\documentclass[11pt]{amsart}
\usepackage{amsmath}
\usepackage{amssymb}
\usepackage{amsfonts}
\usepackage{amsthm}
\usepackage{enumerate}
\usepackage[mathscr]{eucal}
\usepackage{verbatim}
\usepackage{amsthm}
\usepackage{amscd}
\usepackage[mathscr]{eucal}
\usepackage{appendix}
\usepackage{tikz}

\numberwithin{equation}{section}

\newtheorem{innercustomgeneric}{\customgenericname}
\providecommand{\customgenericname}{}

\newcommand{\newcustomtheorem}[2]{\newenvironment{#1}[1]
  {\renewcommand\customgenericname{#2}
   \renewcommand\theinnercustomgeneric{##1}\innercustomgeneric}{\endinnercustomgeneric}}

\newcustomtheorem{customthm}{Theorem}

\newcommand{\newcustomlemma}[2]{\newenvironment{#1}[1]
  {\renewcommand\customgenericname{#2}
   \renewcommand\theinnercustomgeneric{##1} \innercustomgeneric}{\endinnercustomgeneric}}

\newcustomlemma{customlemma}{Lemma}

\newcustomlemma{customproposition}{Proposition}

\newcustomlemma{customclaim}{Claim}

\theoremstyle{plain}
\newtheorem{theorem}{Theorem}[section]
\newtheorem{lemma}[theorem]{Lemma}
\newtheorem{corollary}[theorem]{Corollary}
\newtheorem{proposition}[theorem]{Proposition}

\newtheorem*{remark}{Remark}

\theoremstyle{definition}

\newcommand{\q}{\quad}
\newcommand{\qq}{\qquad}

\newcommand{\bbz}{\mathbb{Z}}

\newcommand{\bbr}{\mathbb{R}}

\newcommand{\bbrn}{\mathbb R^n}

\newcommand{\bbn}{\mathbb{N}}

\newcommand{\II}{\mathcal{I}}

\def\000{\vec{\boldsymbol{0}}}

\newcommand{\supp}{\mathrm{supp}}

\newcommand{\wh}{\widehat}
\newcommand{\wt}{\widetilde}

\newcounter{question}

\newcommand{\bpf}{\begin{proof}}

\newcommand{\epf}{\end{proof}}

\allowdisplaybreaks

\makeatletter
\@namedef{subjclassname@2020}{\textup{2020} Mathematics Subject Classification}
\makeatother

\makeindex         

\begin{document}

\author{Bae Jun Park}
\address{B. Park, Department of Mathematics, Sungkyunkwan University, Suwon 16419, Republic of Korea}
\email{bpark43@skku.edu}

\thanks{The author is supported by POSCO Science Fellowship of POSCO TJ Park Foundation.}

\title{Vector-valued estimates for shifted operators}

\begin{abstract} 
 Shifted variants of (dyadic) Hardy-Littlewood maximal function and Stein's square function have played a significant role in the study of many important operators such as Calder\'on commutators, (bilinear) Hilbert transforms, multilinear multipliers, and multilinear rough singular integrals (e.g. \cite{Do_Park_Sl, Do_Sl2023, Ga_Lie2020, Guo_Hic_Lie_Roo2017, Guo_Roos_Se_Yung2020, LHHLPPY,  Lee_Park2024, Lie2018, Mu_Sc2013}). 
 Estimates for such shifted operators have a certain logarithmic growth in terms of the shift factor, but the optimality of the logarithmic growth has not yet been fully resolved.
  In this article, we provide sharp vector-valued shifted maximal inequality for generalized Peetre's maximal function, from which improved estimates for the above shifted operators follow with optimal logarithmic growths in a new way. We also obtain a vector-valued maximal inequality for the shifted (dyadic) Hardy-Littlewood maximal operator.
\end{abstract}

\maketitle

 \section{Introduction}
 
We denote by $\mathscr{S}(\bbrn)$ the Schwartz class on $\bbrn$ and by $\mathscr{S}'(\bbrn)$ the space of tempered distributions.
For $f\in\mathscr{S}(\bbrn)$ we define its Fourier transform $\wh{f}(\xi):=\int_{\bbrn} f(x) e^{2\pi i\langle x,\xi\rangle}\, dx$ and the inverse Fourier transform 
$f^{\vee}(\xi):=\wh{f}(-\xi)$, and also extend these transforms to the space of tempered distributions.
For $r>0$ let $\mathcal{E}(r)$ denote the space of all tempered distributions whose Fourier transform is supported in the ball of radius $2r$. That is,
$$\mathcal{E}(r):=\big\{f\in\mathscr{S}'(\bbrn): \supp(\wh{f})\subset \{\xi\in\bbrn: |\xi|<2r\} \big\}.$$
 We recall that for each $f\in\mathcal{E}(r)$
 $$f=f\ast \Gamma_r \q \text{ in the sense of tempered distribution} $$
 for a Schwartz function $\Gamma$ whose Fourier transform is equal to $1$ on the ball of radius $2r$, centered at $0$ and is supported in a larger ball. Since a convolution between a tempered distribution and a Schwartz function is a smooth function, we may regard $f\in\mathcal{E}(r)$ as a smooth function with $\supp(\wh{f})\subset \{\xi\in\bbrn: |\xi|<2r\}$.

For a locally integrable function $f$ on $\bbrn$, we define  the classical Hardy-Littlewood maximal function of $f$ by
 \begin{equation*}
 \mathcal{M}f(x):=\sup_{Q: x\in Q}{\frac{1}{|Q|}\int_Q{|f(y)|}dy}
 \end{equation*} 
 where the supremum is taken over all cubes in $\bbrn$ with sides parallel to the axes containing the point $x$. For $0<t<\infty$, we also define its $L^t$ version by
 $\mathcal{M}_tf:=\big(  \mathcal{M}(|f|^t) \big)^{1/t}$ for locally $t$-th power integrable function $f$ on $\bbrn$.
Then the Fefferman-Stein vector-valued maximal inequality \cite{Fe_St1971} states that for $0<t<p,q<\infty$,
\begin{equation}\label{hlmax}
\Big\Vert  \Big(\sum_{k}{(\mathcal{M}_{t}f_k)^q}\Big)^{1/{q}} \Big\Vert_{L^p(\bbrn)} \lesssim  \Big\Vert \Big( \sum_{k}{|f_k|^q}  \Big)^{1/{q}}  \Big\Vert_{L^p(\bbrn)}.
\end{equation}  
The inequality (\ref{hlmax}) also holds for $0<p\leq \infty$ and $q=\infty$.
Here and in the sequel, the symbol $A\lesssim B$ indicates that $A\le CB$ for some constant $C>0$ independent of the variable quantities $A$ and $B$, and $A\sim B$ if $A\lesssim B$ and $B\lesssim A$ hold simultaneously.
For $k\in\mathbb{Z}$ and $\sigma>0$ we recall  Peetre's maximal operator $\mathfrak{M}_{\sigma,2^k}$ defined by the formula
\begin{equation*}
\mathfrak{M}_{\sigma,2^k}f(x):=\sup_{z\in\bbrn}{\frac{|f(x-z)|}{(1+2^k|z|)^{\sigma}}}=\sup_{z\in\bbrn}{\frac{|f(z)|}{(1+2^k|x-z|)^{\sigma}}}.
\end{equation*}
It is known in \cite{Pe1975} that for $f\in \mathcal{E}(A2^k)$,
\begin{equation*}
\mathfrak{M}_{n/t,2^k}f(x)\lesssim_A \mathcal{M}_tf(x) \qquad \text{ uniformly in }~ k
\end{equation*} 
and accordingly, if $\sigma>n/\min{(p,q)}$ and $f_k\in\mathcal{E}(A2^k)$ for some $A>0$, then the maximal inequality
\begin{equation}\label{maximal1}
\Big\Vert \Big\{\mathfrak{M}_{\sigma,2^k}f_k\Big\}_{k\in\mathbb{Z}}\Big\Vert_{L^p(\ell^q)}\lesssim_A \big\Vert \{ f_k\}_{k\in\mathbb{Z}}\big\Vert_{L^p(\ell^q)}
\end{equation}
holds. 

The author \cite{ParkIUMJ} introduced a generalization of Peetre's maximal function 
$$\mathfrak{M}_{\sigma,2^k}^tf(x):= 2^{kn/t}\bigg\Vert \frac{f(x-\cdot)}{(1+2^k|\cdot|)^{\sigma}}\bigg\Vert_{L^t(\bbrn)}, \q 0<t<\infty$$
so that the classical maximal operator can be interpreted as $\mathfrak{M}_{\sigma,2^k}f=\mathfrak{M}^{\infty}_{\sigma,2^k}f.$
We note that for each $\sigma>0$ and $k\in\bbz$, 
\begin{equation}\label{decreasingprom}
\mathfrak{M}_{\sigma,2^k}^tf(x)\lesssim_A \mathfrak{M}_{\sigma,2^k}^sf(x), \q \text{ uniformly in }~k,
\end{equation}
if $0<s\le t\le \infty$ and $f\in\mathcal{E}(A2^k)$ for $A>0$.
Furthermore, for $\sigma>n/\min{(p,q,t)}$,
 $$\mathfrak{M}_{\sigma,2^k}^tf(x)\lesssim_A \mathcal{M}_tf(x) \q \text{ uniformly in }~k,$$
and consequently, \eqref{maximal1} can be strengthened as
\begin{equation}\label{maximal2}
\Big\Vert \Big\{\mathfrak{M}_{\sigma,2^k}^tf_k\Big\}_{k\in\mathbb{Z}}\Big\Vert_{L^p(\ell^q)}\lesssim_A \big\Vert \{ f_k\}_{k\in\mathbb{Z}}\big\Vert_{L^p(\ell^q)},
\end{equation} provided that $\sigma>n/\min{(p,q,t)}$.
In \cite{Park2019, ParkIUMJ}, the above vector-valued maximal inequality has been extended to the case $p=\infty$ and $0<q<\infty$ in the scale of Triebel-Lizorkin spaces.
 Suppose that $0<q<\infty$, $0<t\le \infty$, $\sigma>n/\min{(q,t)}$, and $f_k\in\mathcal{E}(A2^k)$ for some $A>0$. Then we have
\begin{equation}\label{maximal2inf}
\sup_{P\in\mathcal{D}}{\bigg( \frac{1}{|P|}\int_P{\sum_{k=-\log_2{\ell(P)}}^{\infty}{\Big(\mathfrak{M}^t_{\sigma,2^k}f_k(x) \Big)^q}}dx\bigg)^{1/q}}\lesssim\sup_{P\in\mathcal{D}}{\bigg( \frac{1}{|P|}\int_P{\sum_{k=-\log_2{\ell(P)}}^{\infty}{\big|f_k(x) \big|^q}}dx\bigg)^{1/q}}
\end{equation} where the supremums range over all dyadic cubes in $\bbrn$ and $\ell(P)$ denotes the side-length of the cube $P$.
We remark that (\ref{maximal2inf}) does not hold when $\mathfrak{M}^t_{\sigma,2^k}f_k$ is replaced by $\mathcal{M}_rf_k$ for all $0<r<\infty$.
As an application of (\ref{maximal2inf}),  we also obtain
\begin{align}\label{embeddinginfty}
\big\Vert \big\{f_k\big\}_{k\in\bbz}\big\Vert_{L^{\infty}(\ell^{\infty})}&\lesssim \sup_{P\in\mathcal{D}}{\bigg( \frac{1}{|P|}\int_P{\sum_{k=-\log_2{\ell(P)}}^{\infty}{|f_k(x) |^{q}}}dx\bigg)^{1/q}}.
\end{align}
See \cite{Park2019, ParkIUMJ} for more details.

\hfill

Now we consider shifted variants of maximal operators, which are the main objects of this paper.
Let $y\in\bbrn$ be an arbitrary point, which will serve as a shift factor of several operators.
For $0<t<\infty$, we define
$$\mathfrak{M}_{\sigma,2^k,y}^{t}f(x):=\mathfrak{M}_{\sigma,2^k}^{t}f(x-2^{-k}y)=2^{kn/t}\bigg\Vert \frac{\big|f(x-\cdot) \big|}{\big( 1+|2^k\cdot -y|\big)^{\sigma}}\bigg\Vert_{L^t(\bbrn)}$$
and
$$\mathfrak{M}^{\infty}_{\sigma,2^k,y}f(x)=\mathfrak{M}_{\sigma,2^k,y}f(x):=\mathfrak{M}_{\sigma,2^k}f(x-2^{-k}y)=\sup_{z\in\bbrn} \frac{\big|f(x-z) \big|}{\big( 1+|2^kz -y|\big)^{\sigma}}.$$
Appealing to \eqref{decreasingprom}, 
if $0<s\le t\le \infty$ and $f\in\mathcal{E}(A2^k)$ for some $A>0$, then
we have
\begin{equation}\label{decreasingshiftn}
\mathfrak{M}_{\sigma,2^k,y}^tf(x)\lesssim_A \mathfrak{M}_{\sigma,2^k,y}^sf(x) \q \text{ uniformly in ~ $k$}
\end{equation}
where the constant in the inequality is independent of $y\in\bbrn$.

\hfill

The main result of this paper is the following sharp vector-valued maximal inequalities, which generalize \eqref{maximal2} and \eqref{maximal2inf} to  the shifted maximal operator $\mathfrak{M}_{\sigma,2^k,y}^{t}$.
\begin{theorem}\label{maincorol}
Let $0<p,q,t\le \infty$, $y\in\bbrn$, and $\sigma> n/\min{(p,q,t)}$.
For each $k\in\bbn$, let $f_k\in \mathcal{E}(A2^k)$ for some $A>0$.
\begin{enumerate}
\item When $0<p<\infty$ or $p=q=\infty$, 
\begin{equation}\label{rcondition}
r\ge |1/p-1/q|
\end{equation}
if and only if
\begin{equation}\label{finitepest}
\Big\Vert \Big\{\mathfrak{M}_{\sigma,2^k,y}^{t}f_k \Big\}_{k\in\bbz}\Big\Vert_{L^p(\ell^q)}\lesssim \big(\ln{(e+|y|)}\big)^{r}\big\Vert \big\{ f_k\big\}_{k\in\bbz}\big\Vert_{L^p(\ell^q)}
\end{equation}
\item When $p=\infty$ and $0<q<\infty$, 
\begin{equation}\label{rcondition11}
r\ge 1/q
\end{equation}
if and only if
\begin{align}\label{pinfqfiniteesto}
&\sup_{P\in\mathcal{D}}\bigg( \frac{1}{|P|}\int_P \sum_{k:2^k\ell(P)\ge 1}\Big(\mathfrak{M}_{\sigma,2^k,y}^{t}f_k(x) \Big)^q\;dx \bigg)^{1/q}\nonumber\\
&\lesssim \big( \ln{(e+|y|)}\big)^{r} \sup_{P\in\mathcal{D}}\bigg( \frac{1}{|P|}\int_P \sum_{k:2^k\ell(P)\ge 1}\big|f_k(x) \big|^q\;dx \bigg)^{1/q}.
\end{align}

\end{enumerate}
\end{theorem}

In order to prove one direction (``only if\," part) of Theorem \ref{maincorol}, we will actually prove the following proposition.
\begin{proposition}\label{maintheorem}
Let $0<p,q,t \le \infty$ and $y\in\bbrn$.
Suppose that
$$\begin{cases}
t\le q & \text{if} \q q\le p\\
t<p  & \text{if} \q p<q
\end{cases} \qq \text{ and }\qq \sigma>n/t.$$
\begin{enumerate}
\item If $0<p<\infty$ or $p=q=\infty$, then we have
\begin{align*}
\Big\Vert \Big\{\mathfrak{M}_{\sigma,2^k,y}^{t}f_k \Big\}_{k\in\bbz}\Big\Vert_{L^p(\ell^q)}\lesssim \big(\ln{(e+|y|)}\big)^{|1/p-1/q|}\big\Vert \big\{ f_k\big\}_{k\in\bbz}\big\Vert_{L^p(\ell^q)}
\end{align*}
\item If $p=\infty$ and $0<q<\infty$, then we have
\begin{align}\label{pinfqfiniteest}
&\sup_{P\in\mathcal{D}}\bigg( \frac{1}{|P|}\int_P \sum_{k:2^k\ell(P)\ge 1}\Big(\mathfrak{M}_{\sigma,2^k,y}^{t}f_k(x) \Big)^q\;dx \bigg)^{1/q}\nonumber\\
&\lesssim \big( \ln{(e+|y|)}\big)^{1/q}\bigg(\big\Vert \{f_k\}_{k\in\bbz}\big\Vert_{L^{\infty}(\ell^{\infty})}+\sup_{P\in\mathcal{D}}\bigg( \frac{1}{|P|}\int_P \sum_{k:2^k\ell(P)\ge 1}\big|f_k(x) \big|^q\;dx \bigg)^{1/q} \bigg).
\end{align}

\end{enumerate}
\end{proposition}
We point out that the Fourier support condition of $f_k$ is unnecessary in the proposition, unlike Theorem \ref{maincorol}.
Once Proposition \ref{maintheorem} is established, one direction of Theorem \ref{maincorol}    follows immediately from  \eqref{decreasingprom} and \eqref{embeddinginfty} assuming $f_k\in\mathcal{E}(A2^k)$ for each $k\in\bbz$.

For the opposite direction (``if\," part in Theorem \ref{maincorol}), we will construct counterexamples to verify the validity of the conditions \eqref{rcondition} and \eqref{rcondition11} for the estimates \eqref{finitepest} and \eqref{pinfqfiniteesto} to hold, respectively. This will be contained in Section \ref{counterexamplesec}.

\hfill

\subsection{Shifted dyadic maximal function estimates}
For any $y\in\mathbb{R}^n$ and for any dyadic cubes $Q\in\mathcal{D}$, we use the notation
 $$Q(y):=Q+\ell(Q)y$$
and then define the shifted (dyadic) Hardy-Littlewood maximal operator $\mathscr{M}^y$ by
$$\mathscr{M}^yf(x):=\sup_{Q\in\mathcal{D}:x\in Q}\frac{1}{|Q|}\int_{Q(y)} |f(z)|\; dz$$
and for $0<t<\infty$
$$\mathscr{M}_t^yf(x):=\Big(\mathscr{M}^y\big(|f|^t \big)(x)\Big)^{1/t}.$$

As an application of Proposition \ref{maintheorem}, we obtain the following shifted maximal inequality.
\begin{theorem}\label{dyadmaxtheorem}
Let $0<p<\infty$, $y\in\bbrn$, and $0<t<p$.
Then we have
$$\big\Vert \mathscr{M}^y_tf\big\Vert_{L^p(\bbrn)}\lesssim \big( \ln{(e+|y|)}\big)^{1/p}\Vert f\Vert_{L^p(\bbrn)}.$$
\end{theorem}
Analogous estimates in one dimensional space $\bbr$ were already studied in the paper of Muscalu  \cite[Theorem 4.1]{Mu2014} and in the book of Stein \cite[Chapter II, 5.10]{St1993} based on the Calder\'on-Zygmund theory, and general $n$-dimensional extensions of such estimates are introduced in the book of Grafakos \cite[Proposition 7.5.1]{Gr2014} and in the paper of Lee et al. \cite[(3.2)]{LHHLPPY}, but they do not pursue the optimal growth in terms of $y$ as their results require the growth
$(\ln{(e+|y|)})^{n/p}$ or $\ln{(e+|y|)}$ for $p>1$ and $t=1$. 
Theorem \ref{dyadmaxtheorem} improves the previous results, providing the optimal logarithmical bound.
This sharpness will follow from that of Theorem \ref{phikyastfest} below.
We would like to emphasize that our approach to proving Theorem \ref{dyadmaxtheorem}, mainly relying on Proposition \ref{maintheorem} without interpolation techniques, is 
completely different from those of the earlier works.

\hfill

We also prove a vector-valued version of the above maximal inequality.
\begin{theorem}\label{vectorvalueddyadic}
Let $0<p<\infty$, $0<q\le \infty$,  $y\in\bbrn$, and $0<t<p,q$. 
Then we have
\begin{equation}\label{vecdyashiftest}
\big\Vert \big\{ \mathscr{M}_t^yf_k\big\}_{k\in\bbz}\big\Vert_{L^p(\ell^q)}\lesssim \big(\ln{(e+|y|)} \big)^{\tau_{p,q}}\Vert \{ f_k\}_{k\in\bbz}\Vert_{L^p(\ell^q)}
\end{equation}
where 
$\tau_{p,q}$ is equal to $\frac{1}{p}$ if $p\le q$, but $\tau_{p,q}$ is any number greater than $\frac{1}{q}$ if $q<p$.
\end{theorem}
Clearly, the above inequality holds for $p=q=\infty$ as well. 
A weaker version of such a maximal inequality  appeared in the paper of Guo, Hickman, Lie, and Roos \cite[Theorem 3.1]{Guo_Hic_Lie_Roo2017} who obtained the $L^p(\ell^q)$ estimates for $\mathscr{M}^y$ with the logarithmic exponent $2$,  via a remarkable weighted norm inequality, in the one dimensional setting. However, as mentioned in the paper, they did not consider the best logarithmic exponent as the number $2$ is enough for their purpose.

The proof of Theorem \ref{vectorvalueddyadic} is based on the Calder\'on-Zygmund theory for vector-valued function spaces with Theorem \ref{dyadmaxtheorem} being a kind of initial estimate, but a certain dilation trick is necessary when $q<p$. 

\begin{remark}
$\tau_{p,q}=\frac{1}{p}$ if $p\le q$ and it is the best logarithmic exponent for \eqref{vecdyashiftest} to hold, considering the sharpness of Theorem \ref{dyadmaxtheorem}.  When $q<p$, $\tau_{p,q}$ could be arbitrary close to $1/q$ in Theorem \ref{vectorvalueddyadic}. One might guess that the optimal logarithmic exponent would be $\frac{1}{\min{(p,q)}}$ for the vector-valued maximal inequality \eqref{vecdyashiftest},
 but we could not verify it here.
\end{remark}

\hfill

\subsection{Shifted function estimates for Hardy spaces and $BMO$}

Let $\phi$ be a  Schwartz function on $\bbrn$ satisfying 
$$\supp(\wh{\phi})\subset \{\xi\in\bbrn: |\xi|\le 2\}, \q \wh{\phi}(\xi)=1~\text{ for }~|\xi|\le 1.$$
Let $\psi:=\phi-2^{-n}\phi(2^{-1}\cdot)$ so that
$$\supp(\wh{\psi})\subset \big\{\xi\in\bbrn: 2^{-1}\le |\xi|\le 2 \big\},\q\text{ and }\q
\sum_{k\in\bbz}\wh{\psi_k}(\xi)=1 ~ \text{ for } ~ \xi\not= 0$$
where we set $\phi_k:=2^{kn}\phi(2^k\cdot)$ and $\psi_k:=2^{kn}\psi(2^k\cdot)$ for each $k\in\bbz$.
Then the (real) Hardy space $H^p(\bbrn)$ ($0<p<\infty$), consists of tempered distributions $f$ on $\bbrn$ such that 
\begin{equation*}
\Vert f \Vert_{H^p(\bbrn)}:= \Big\Vert \sup_{k\in\bbz} \big| \phi_k\ast f \big|\Big\Vert_{L^p(\bbrn)}
\end{equation*}
is finite.
The space $H^p$ is same as the Lebesgue space $L^p$ for $1<p<\infty$, and naturally extends the interpolation scale of the $L^p$ spaces, $1<p<\infty$, to all $0<p<\infty$.
As a consequence of Theorem \ref{maincorol}, we obtain the following shifted estimate for Hardy space.
For each $k\in \bbz$ and $y\in\bbrn$ let $$(\phi_k)^y:=\phi_k(\cdot-2^{-k}y)=2^{kn}\phi(2^k\cdot-y).$$
Then the following theorem is immediately obtained from Theorem \ref{maincorol} if we notice the pointwise estimate
\begin{equation}\label{phikyastfpoint}
\big| (\phi_k)^y\ast f(x)\big|=\big| \phi_k\ast f(x-2^{-k}y)\big|\le \mathfrak{M}_{\sigma,2^k,y}\big(\phi_k\ast f\big)(x)
\end{equation} for any $\sigma>0$.
\begin{theorem}\label{phikyastfest}
Let $0<p< \infty$ and $y\in \bbrn$.
Then we have
$$\bigg\Vert \sup_{k\in\bbz}\big| (\phi_k)^y\ast f\big|\bigg\Vert_{L^p(\bbrn)}\lesssim \big( \ln{(e+|y|)}\big)^{1/p}\Vert f\Vert_{H^p(\bbrn)}.$$
Moreover,  the number  $1/p$ is the best logarithmic exponent for which the inequality holds. 
\end{theorem}
The sharpness of the theorem will be discussed in Section 6.

\hfill

The Hardy space can be characterized in terms of Stein's square function as
\begin{equation}\label{hardycha}
\Vert f \Vert_{H^p(\bbrn)}\sim \bigg\Vert \bigg( \sum_{k\in\bbz}\big| \psi_k\ast f \big|^2\bigg)^{1/2}\bigg\Vert_{L^p(\bbrn)},
\end{equation}
which is known as Littlewood-Paley theory for Hardy space. 
As an end-point case $p=\infty$ of the Lebesgue space $L^p$ (or Hardy space $H^p$), we typically consider $BMO$ space, which is, in particular, the dual of $H^1$. The $BMO$ space can be also characterized in terms of the square function, by
\begin{equation}\label{bmocha}
\Vert f\Vert_{BMO}\sim \sup_{P\in\mathcal{D}}\bigg(\frac{1}{|P|}\int_P \sum_{k:2^k\ell(P)\ge 1}\big| \psi_k\ast f(x)\big|^2\;dx \bigg)^{1/2},
\end{equation}
which originated from Carleson measure characterization.
For each $k\in\bbz$ and $y\in\bbrn$, we now define
$$(\psi_k)^y:=\psi_k(\cdot-2^{-k}y)=2^{kn}\psi(2^k\cdot-y).$$
As an additional consequence of Theorem \ref{maincorol}, we obtain the following vector-valued estimates for shifted Littlewood-Paley operators.
\begin{theorem}\label{shiftedvectorest}
Let $0<p,q\le \infty$ and $y\in\bbrn$.
\begin{enumerate}
\item If $0<p<\infty$ or $p=q=\infty$, then we have
\begin{equation}\label{psikyastflplqest}
\big\Vert \big\{ (\psi_k)^y\ast f\big\}_{k\in\bbz}\big\Vert_{L^p(\ell^q)}\lesssim \big( \ln{(e+|y|)}\big)^{|1/p-1/q|}\big\Vert \big\{ \psi_k\ast f\big\}_{k\in\bbz}\big\Vert_{L^p(\ell^q)}.
\end{equation}
\item If $p=\infty$ and $0<q<\infty$, 
then we have
\begin{align}\label{inftypsikyastflplqest}
&\sup_{P\in\mathcal{D}}\bigg( \frac{1}{|P|}\int_P \sum_{k:2^k\ell(P)\ge 1}\big|(\psi_k)^y\ast f(x) \big|^q\; dx\bigg)^{1/q}\nonumber\\
&\lesssim \big( \ln{(e+|y|)}\big)^{1/q}\sup_{P\in\mathcal{D}}\bigg( \frac{1}{|P|}\int_P \sum_{k:2^k\ell(P)\ge 1}\big|\psi_k \ast f(x) \big|^q\; dx\bigg)^{1/q}.
\end{align}
\end{enumerate}
Moreover, the numbers $|1/p-1/q|$ and $1/q$ are the best logarithmic exponents for which the inequalities hold.
\end{theorem}
Similar to \eqref{phikyastfpoint}, we have
$$\big| (\psi_k)^y\ast f(x)\big| \le \mathfrak{M}_{\sigma,2^k,y}\big(\psi_k\ast f\big)(x)$$
and then Theorem \ref{maincorol} leads to the inequalities in Theorem \ref{shiftedvectorest}.
The sharpness can be achieved by analogous (or essentially same) constructions of counterexamples used in the proof of ``if\," part of Theorem \ref{maincorol}. We will include this in Section \ref{counterexamplesec}.

\hfill

Now, characterizations \eqref{hardycha} and \eqref{bmocha} yield the following sharp square function estimates for Hardy space and $BMO$ space, which generalize one direction of the characterizations themselves.
\begin{corollary}\label{shiftedsquare}
Let $0<p\le \infty$ and $y\in\bbrn$.
\begin{enumerate}
\item If $0<p<\infty$, then
$$\bigg\Vert \bigg( \sum_{k\in\bbz}\big| (\psi_k)^y\ast f\big|^2\bigg)^{1/2}\bigg\Vert_{L^p(\bbrn)}\lesssim \big( \ln{(e+|y|)}\big)^{|1/p-1/2|}\Vert f\Vert_{H^p(\bbrn)}.$$
\item If $p=\infty$, then
$$\sup_{P\in\mathcal{D}}\bigg( \frac{1}{|P|}\int_P \sum_{k:2^k\ell(P)\ge 1}\big|(\psi_k)^y\ast f(x) \big|^2\; dx\bigg)^{1/2}\lesssim \big( \ln{(e+|y|)}\big)^{1/2}\Vert f\Vert_{BMO}.$$
\end{enumerate}
\end{corollary}

The optimality of the logarithmic exponents $|1/p-1/2|$ and $1/2$ in the estimates is guaranteed by the sharpness of Theorem \ref{shiftedvectorest}.
When $1<p<\infty$, a discrete version of such a shifted square function estimate appeared in \cite[Theorem 4.6]{Mu_Sc2013} for one-dimensional case and higher-dimensional ones are also studied in  \cite[Corollary 1]{Gr_Oh2014}, but both are valid with the logarithmic exponent $|1/p-1/2|$ replaced by $1$. 
See also \cite[Lemma 3]{Lie2018}.
 Recently, these results have  been improved in \cite[Lemma 4]{Do_Park_Sl} to have logarithmic exponent $|1/p-1/2|$ and extended to (endpoint) Banach range $1\le p\le \infty$,  using certain interpolation and duality arguments, which cannot be applied to the case $0<p<1$.
  In conclusion, 
Corollary \ref{shiftedsquare} provides the final form of all of the earlier estimates with optimal logarithmic bounds in the full range $0<p\le \infty$.


\hfill

{\bf Organization.} In Section \ref{pfpropo12} we provide the proof of Proposition \ref{maintheorem}.
The proof of Theorem \ref{dyadmaxtheorem} is given in Section \ref{proofthm13}.
Section \ref{keylemmathm14} contains some lemmas which will be very useful tools in dealing with vector-valued inequalities. By using the lemmas, we prove Theorem \ref{vectorvalueddyadic} in Section \ref{proofthm14}.
  Section \ref{counterexamplesec} is devoted to constructing several counterexamples which verify the sharpness of Theorems \ref{maincorol}, \ref{phikyastfest}, and \ref{shiftedvectorest}.

\hfill

\section{Proof of Proposition \ref{maintheorem}}\label{pfpropo12}
Assume $\sigma>n/t$.
We first observe that for all $0< t\le r< \infty$,
\begin{align}\label{basiclpest}
\Big\Vert \mathfrak{M}_{\sigma,2^k,y}^{t}f_k \Big\Vert_{L^r(\bbrn)}&\le 2^{kn/t}\bigg(\int_{\bbrn}  \Big(\int_{\bbrn}\frac{|f_k(x-z)|^r}{(1+|2^k z-y|)^{\sigma r}} \, dx\Big)^{t/r}     \, dz \bigg)^{1/t}\nonumber\\
&= \bigg(\int_{\bbrn}\frac{2^{kn}}{(1+2^k|z|)^{\sigma t}}\Big(\int_{\bbrn}\big| f_k(x-z)\big|^r\; dx  \Big)^{t/r}\; dz   \bigg)^{1/t}  \sim \Vert f_k\Vert_{L^r(\bbrn)}
\end{align}
where the inequality follows from Minkowski's inequality if $r>t$.
 Consequently, we have
 \begin{align*}
 \Big\Vert \Big\{  \mathfrak{M}_{\sigma,2^k,y}^tf_k \Big\}_{k\in\bbz}\Big\Vert_{L^p(\ell^p)}&= \bigg\Vert \bigg\{\Big\Vert \mathfrak{M}_{\sigma,2^k,y}^{t}f_k\Big\Vert_{L^p(\bbrn)} \bigg\}_{k\in\bbz}\bigg\Vert_{\ell^p}\\
 &\lesssim \Big\Vert \Big\{\big\Vert f_k\big\Vert_{L^p(\bbrn)} \Big\}_{k\in\bbz}\Big\Vert_{\ell^p}=\big\Vert \big\{ f_k\big\}_{k\in\bbz}\big\Vert_{L^p(\ell^p)}
 \end{align*} for $0<p\le \infty$ and $t\le p$. Here, the case $p=\infty$ is trivial.
 Therefore, we will consider only the cases when $p\not= q$.

 \subsection{The case when $0<q<p< \infty$}\label{0<q<p<inf}
 In this case, the proof mainly depends on a fundamental inequality of Fefferman and Stein  \cite{Fe_St1972} for sharp maximal function.
For a locally integrable function $f$ on $\bbrn$, we define the sharp (dyadic) maximal function $\mathscr{M}^{\sharp}f$ by
 $$\mathscr{M}^{\sharp}f(x):=\sup_{Q\in \mathcal{D}:x\in Q}\inf_{c_Q\in \mathbb{C}}\bigg(\frac{1}{|Q|}\int_Q \big| f(y)-c_Q\big|\; dy\bigg)$$
 where the supremum is taken over all dyadic cubes in $\bbrn$ containing $x$.
Then for $1\le p_0\le p<\infty$ and $f\in L^{p_0}(\bbrn)$ we have
\begin{equation}\label{flpsharpmax}
\Vert f\Vert_{L^p(\bbrn)}\lesssim \big\Vert \mathscr{M}^{\sharp}f\big\Vert_{L^p(\bbrn)}.
\end{equation} 
In fact, non-dyadic sharp maximal function was used in the original version of \eqref{flpsharpmax}, but as mentioned in \cite[Section 2.5]{Park_BLMS}, \eqref{flpsharpmax} also holds with dyadic sharp maximal one.

We may assume that $t=q$ and $\sigma>n/q$.
 We first write
 \begin{align*}
  \Big\Vert \Big\{\mathfrak{M}_{\sigma,2^k,y}^{q}f_k \Big\}_{k\in\bbz}\Big\Vert_{L^p(\ell^q)}= \bigg\Vert \sum_{k\in\bbz}\Big( \mathfrak{M}_{\sigma,2^k,y}^qf_k \Big)^q\bigg\Vert_{L^{p/q}(\bbrn)}^{1/q}
 \end{align*} and then
 this is controlled, via \eqref{flpsharpmax}, by a constant times
$$ \bigg\Vert \mathscr{M}^{\sharp}\bigg( \sum_{k\in\bbz}\Big( \mathfrak{M}_{\sigma,2^k,y}^qf_k \Big)^q\bigg)\bigg\Vert_{L^{p/q}(\bbrn)}^{1/q}.$$
 The sharp maximal function inside the $L^{p/q}$ norm is bounded by 
 $$2\bigg( \mathscr{M}_q^{\sharp,(1)}\Big(\Big\{ \mathfrak{M}_{\sigma,2^k,y}^qf_k \Big\}_{k\in\bbz} \Big) (x)\bigg)^q+\bigg(\mathscr{M}_q^{\sharp,(2)}\Big(\Big\{ \mathfrak{M}_{\sigma,2^k,y}^qf_k \Big\}_{k\in\bbz} \Big) (x)\bigg)^q$$
 where
 \begin{equation}\label{mqs,1se}
 \mathscr{M}_q^{\sharp,(1)}\Big(\Big\{ \mathfrak{M}_{\sigma,2^k,y}^qf_k \Big\}_{k\in\bbz} \Big) (x):=\sup_{P\in\mathcal{D}:x\in P}\bigg( \frac{1}{|P|}\int_P \sum_{k:2^k\ell(P)\ge 1} \Big( \mathfrak{M}_{\sigma,2^k,y}^qf_k (z) \Big)^q  \; dz\bigg)^{1/q} 
 \end{equation}
  and
  \begin{align*}
&\mathscr{M}_q^{\sharp,(2)}\Big(\Big\{ \mathfrak{M}_{\sigma,2^k,y}^qf_k \Big\}_{k\in\bbz} \Big) (x)\\
&:=  \sup_{P\in\mathcal{D}:x\in P}\bigg( \frac{1}{|P|}\int_P\frac{1}{|P|}\int_P \sum_{k:2^k\ell(P)<1}\Big| \Big( \mathfrak{M}_{\sigma,2^k,y}^qf_k (z)\Big)^q-\Big( \mathfrak{M}_{\sigma,2^k,y}^qf_k (u)  \Big)^q       \Big| \; dz \; du \bigg)^{1/q}.
  \end{align*}
 Then we will actually prove that
  \begin{equation}\label{qpclaim1}
\bigg\Vert \mathscr{M}_q^{\sharp,(1)}\Big(\Big\{ \mathfrak{M}_{\sigma,2^k,y}^qf_k \Big\}_{k\in\bbz} \Big)\bigg\Vert_{L^p(\bbrn)}\lesssim \big( \ln{(e+|y|)}\big)^{1/q-1/p}\big\Vert \big\{ f_k\big\}_{k\in\bbz}\big\Vert_{L^p(\ell^q)}
 \end{equation}
 and
  \begin{equation}\label{qpclaim2}
\bigg\Vert \mathscr{M}_q^{\sharp,(2)}\Big(\Big\{ \mathfrak{M}_{\sigma,2^k,y}^qf_k \Big\}_{k\in\bbz} \Big)\bigg\Vert_{L^p(\bbrn)}\lesssim \big\Vert \big\{ f_k\big\}_{k\in\bbz}\big\Vert_{L^p(\ell^q)}.
 \end{equation}

\subsubsection{Proof of \eqref{qpclaim1}}
 We decompose the left-hand side of \eqref{qpclaim1} as
 \begin{align}\label{lpsuppind}
& \bigg\Vert \sup_{P\in\mathcal{D}:x\in P}\bigg(\frac{1}{|P|}\int_P{\sum_{k:2^k\ell(P)\ge 1}{\Big(\mathfrak{M}_{\sigma,2^k,y}^q\big(\chi_{P^{*}} f_k \big)(z) \Big)^q}} \;dz \bigg)^{1/q}\bigg\Vert_{L^p(x)}\nonumber\\
 &\q + \bigg\Vert \sup_{P\in\mathcal{D}:x\in P}\bigg(\frac{1}{|P|}\int_P{\sum_{k:1\le 2^k\ell(P)<2+|y|}{\Big(\mathfrak{M}_{\sigma,2^k,y}^q\big( \chi_{(P^*)^c}f_k\big)(z) \Big)^q}} \;dz \bigg)^{1/q}\bigg\Vert_{L^p(x)}\nonumber\\
  &\q \q+ \bigg\Vert \sup_{P\in\mathcal{D}:x\in P}\bigg(\frac{1}{|P|}\int_P{\sum_{k: 2^k\ell(P)\ge 2+|y|}{\Big(\mathfrak{M}_{\sigma,2^k,y}^q\big( \chi_{(P^*)^c}f_k\big)(z) \Big)^q}} \;dz \bigg)^{1/q}\bigg\Vert_{L^p(x)}\nonumber\\
 &=:\mathrm{I}_1^{(p,q)}+\mathrm{I}_2^{(p,q)}+\mathrm{I}_3^{(p,q)}
 \end{align}
 where $P^*$ denotes the concentric dilate of $P$ with $\ell(P^*)=10\sqrt{n}\ell(P)$.
 The three terms are treated separately below and we point out that
 the logarithmic increasing factor $\big( \ln{(e+|y|)}\big)^{1/q-1/p}$ is needed only for the second one $\mathrm{I}_2^{(p,q)}$.

 For the estimate for $\mathrm{I}_1^{(p,q)}$, we employ the inequality \eqref{basiclpest} to obtain
 \begin{align*}
 \int_P{\sum_{k:2^k\ell(P)\ge 1}{\Big(\mathfrak{M}_{\sigma,2^k,y}^q\big(\chi_{P^{*}} f_k \big)(z) \Big)^q}} \;dz &\le \sum_{k:2^k\ell(P)\ge 1}\Big\Vert \mathfrak{M}_{\sigma,2^k,y}^{q}\big(\chi_{P^*}f_k\big)\Big\Vert_{L^q(\bbrn)}^q\\
 &\lesssim \int_{P^*} \sum_{k:2^k\ell(P)\ge 1}\big| f_k(z)\big|^q\; dz,
 \end{align*}
 which proves that
 \begin{equation}\label{II1pest}
 \mathrm{I}_1^{(p,q)}\lesssim \bigg\Vert \sup_{P\in\mathcal{D}:x\in P}\bigg(\frac{1}{|P|}\int_{P^*}\sum_{k:2^k\ell(P)\ge 1}\big| f_k(z)\big|^q \; dz\bigg)^{1/q}\bigg\Vert_{L^p(x)}\lesssim \Big\Vert   \mathscr{M}_q^{\sharp,(1)}\big( \{f_k\}_{k\in\bbz}\big)   \Big\Vert_{L^p(\bbrn)}
 \end{equation}
 where $\mathscr{M}_q^{\sharp,(1)}\big( \{f_k\}_{k\in\bbz}\big)$ is defined by replacing $\mathfrak{M}_{\sigma,2^k,y}^qf_k$ in \eqref{mqs,1se} by $f_k$.
 Since
 $$\mathscr{M}_q^{\sharp,(1)}\big( \{f_k\}_{k\in\bbz}\big)(x)\le \mathcal{M}_q \big(\big\Vert \{f_k\}_{k\in\bbz}\big\Vert_{\ell^q} \big)(x),$$
 the maximal inequality \eqref{hlmax} finally yields that
 $$\mathrm{I}_1^{(p,q)}\lesssim \big\Vert \{f_k\}_{k\in\bbz}\big\Vert_{L^p(\ell^q)},$$
 as expected.

 For the next one $\mathrm{I}_2^{(p,q)}$, we bound the sum over $k$ with $1\le2^k\ell(P)<2+|y|$  by
 \begin{align*}
 &\bigg(\sum_{k:1\le 2^k\ell(P)<2+|y|} 1 \bigg)^{1-q/p}    \Big\Vert \Big\{ \mathfrak{M}_{\sigma,2^k,y}^{q}f_k(z)\Big\}_{k\in\bbz} \Big\Vert_{\ell^p}^{q}\\
 &\lesssim \bigg( \big(\ln{(e+|y|)}\big)^{1/q-1/p}\Big\Vert \Big\{\mathfrak{M}_{\sigma,2^k,y}^qf_k(z) \Big\}_{k\in\bbz}\Big\Vert_{\ell^p}\bigg)^q
 \end{align*}
by using H\"older's inequality, and thus $\mathrm{I}_2^{(p,q)}$ is less than a constant times
 \begin{equation*}
 \big(\ln{(e+|y|)}\big)^{1/q-1/p} \bigg\Vert \mathcal{M}_q\bigg( \Big\Vert \Big\{\mathfrak{M}_{\sigma,2^k,y}^qf_k \Big\}_{k\in\bbz}\Big\Vert_{\ell^p}\bigg)\bigg\Vert_{L^p(\bbrn)}.
 \end{equation*}
 Then the inequaliteis \eqref{hlmax} and \eqref{basiclpest} yield that 
 \begin{equation*}
 \bigg\Vert \mathcal{M}_q\bigg( \Big\Vert \Big\{\mathfrak{M}_{\sigma,2^k,y}^qf_k \Big\}_{k\in\bbz}\Big\Vert_{\ell^p}\bigg)\bigg\Vert_{L^p(\bbrn)}\lesssim \Big\Vert \Big\{\mathfrak{M}_{\sigma,2^k,y}^qf_k \Big\}_{k\in\bbz}\Big\Vert_{L^p(\ell^p)}\lesssim \big\Vert \{f_k\}_{k\in\bbz}\big\Vert_{L^p(\ell^p)}.
 \end{equation*} 
 This shows that 
 \begin{equation}\label{II2pest}
  \mathrm{I}_2^{(p,q)} \lesssim \big(\ln{(e+|y|)}\big)^{1/q-1/p} \big\Vert \{f_k\}_{k\in\bbz}\big\Vert_{L^p(\ell^p)}
 \end{equation}
 and then the desired result follows from the embedding $\ell^q\hookrightarrow \ell^p$.

 We now handle the last one $\mathrm{I}_3^{(p,q)}$. For this purpose, we observe that
 if $2^k\ell(P)\ge 2+|y|$, then $|2^{-k}y|\le \ell(P)$ so that for $z\in P$ and $u\in (P^*)^c$,
 \begin{equation*}
 |z-u-2^{-k}y|\ge |z-u|-|2^{-k}y|>|z-u|-\ell(P)\ge \frac{|z-u|}{2}\gtrsim \ell(P).
 \end{equation*}
 This implies that for $z\in P$, 
 \begin{align*}
 \mathfrak{M}_{\sigma,2^k,y}^q\big( \chi_{(P^*)^c}f_k\big)(z)&=\bigg(2^{kn} \int_{(P^*)^c}   \frac{|f_k(u)|^q}{(1+2^k|z-u-2^{-k}y|)^{\sigma q}} \; du  \bigg)^{1/q}\nonumber\\
 &\lesssim \big( 2^k\ell(P)\big)^{-\epsilon_0/q}\bigg( 2^{ kn} \int_{\bbrn}  \frac{ |f_k(u)|^q }{\big( 1+2^k|z-u|\big)^{\sigma q-\epsilon_0}}  \; du \bigg)^{1/q}\nonumber\\
 &\lesssim   \big( 2^k\ell(P)\big)^{-\epsilon_0/q}\mathfrak{M}_{\sigma-\epsilon_0/q,2^k}^qf_k(z)
 \end{align*}
 by  choosing $0<\epsilon_0<\sigma q-n$.
Using this, we obtain
 \begin{align}\label{II3pest}
 \mathrm{I}_3^{(p,q)}&\lesssim \bigg\Vert \sup_{P\in\mathcal{D}:x\in P} \bigg(\frac{1}{|P|}\int_P{\Big\Vert \Big\{\mathfrak{M}_{\sigma-\epsilon_0/q,2^k}^qf_k(z)\Big\}_{k\in\bbz}\Big\Vert_{\ell^{\infty}}^q}    \; dz \bigg)^{1/q}\bigg\Vert_{L^p(x)}\nonumber\\
 &\le \bigg\Vert \mathcal{M}_q \bigg(  \Big\Vert \Big\{ \mathfrak{M}_{\sigma-\epsilon_0/q,2^k}^qf_k\Big\}_{k\in\bbz}\Big\Vert_{\ell^{p}}     \bigg)\bigg\Vert_{L^p(\bbrn)}\lesssim \big\Vert \big\{ f_k\big\}_{k\in\bbz}\big\Vert_{L^p(\ell^{p})}
 \end{align}
  where we applied the embedding $\ell^p\hookrightarrow \ell^{\infty}$,  maximal inequality \eqref{hlmax}, and \eqref{basiclpest} with $\sigma-\epsilon_0/q>n/q$.
 Finally, the embedding $\ell^q\hookrightarrow \ell^{p}$ deduces the desired estimate, completing the proof of  \eqref{qpclaim1}.

 \subsubsection{Proof of \eqref{qpclaim2}}
 To verify \eqref{qpclaim2}, we first observe that
 \begin{align*}
& \Big| \Big( \mathfrak{M}_{\sigma,2^k,y}^qf_k (z)\Big)^q-\Big( \mathfrak{M}_{\sigma,2^k,y}^qf_k (u)  \Big)^q       \Big| \\
&\le \int_{\bbrn}  2^{kn}\big| f_k(v-2^{-k}y)\big|^q  \bigg|\frac{1}{\big( 1+2^k|z-v|\big)^{\sigma q}}-\frac{1}{\big( 1+2^k|u-v|\big)^{\sigma q}} \bigg|  \; dv\\
&\lesssim 2^k|u-z| \int_0^1\int_{\bbrn} 2^{kn}\big| f_k(v-2^{-k}y)\big|^q \frac{1}{\big(1+2^k| rz+(1-r)u-v    | \big)^{\sigma q+1}}        \; dv \; dr.
 \end{align*}
 When $u,z\in P$ and $2^k\ell(P)<1$,
 the preceding expression is controlled by
 $$2^k\ell(P)\Big(\mathfrak{M}_{\sigma+1/q,2^k,y}^qf_k(z)\Big)^q$$
 as
 $|u-z|\lesssim \ell(P)$
 and
 $1+2^k|rz+(1-r)u-v|\gtrsim 1+2^k|z-v|.$
 Therefore,
 \begin{align*}
 \mathscr{M}_q^{\sharp,(2)}\Big(\Big\{ \mathfrak{M}_{\sigma,2^k,y}^qf_k \Big\}_{k\in\bbz} \Big) (x)&\lesssim \sup_{P\in\mathcal{D}: x\in P}\bigg( \frac{1}{|P|}\int_P  \sum_{k:2^k\ell(P)<1} \Big(2^k \ell(P)\mathfrak{M}_{\sigma+1/q,2^k,y}^qf_k(z) \Big)^q  \; dz\bigg)^{1/q}\\
 &\lesssim  \sup_{P\in\mathcal{D}: x\in P}\bigg( \frac{1}{|P|}\int_P\Big\Vert \Big\{ \mathfrak{M}_{\sigma+1/q,2^k,y}^qf_k(z)\Big\}_{k\in\bbz} \Big\Vert_{\ell^{\infty}}^q  \; dz\bigg)^{1/q}\\
 &\le \mathcal{M}_q\bigg( \Big\Vert \Big\{\mathfrak{M}_{\sigma+1/q,2^k,y}^qf_k \Big\}_{k\in\bbz}\Big\Vert_{\ell^{\infty}}\bigg)(x)
 \end{align*}
 and then using the inequality \eqref{hlmax}, the embedding $\ell^p\hookrightarrow \ell^{\infty}$, estimate \eqref{basiclpest}, and the embedding $\ell^q\hookrightarrow \ell^p$, we obtain
 \begin{align*}
& \bigg\Vert \mathscr{M}_q^{\sharp,(2)}\Big(\Big\{ \mathfrak{M}_{\sigma,2^k,y}^qf_k \Big\}_{k\in\bbz} \Big)\bigg\Vert_{L^p(\bbrn)}\\
&\lesssim \bigg\Vert \mathcal{M}_q\bigg( \Big\Vert \Big\{\mathfrak{M}_{\sigma+1/q,2^k,y}^qf_k \Big\}_{k\in\bbz}\Big\Vert_{\ell^{\infty}}\bigg)\bigg\Vert_{L^p(\bbrn)} \lesssim \Big\Vert  \Big\{ \mathfrak{M}_{\sigma+1/q,2^k,y}^qf_k \Big\}_{k\in\bbz}\Big\Vert_{L^p(\ell^{\infty})}\\
&\le \bigg( \sum_{k\in\bbz}\Big\Vert    \mathfrak{M}_{\sigma+1/q,2^k,y}^qf_k       \Big\Vert_{L^p(\bbrn)}^p\bigg)^{1/p}\lesssim \big\Vert \{f_k\}_{k\in\bbz}\big\Vert_{L^p(\ell^p)}\lesssim \big\Vert \{f_k\}_{k\in\bbz}\big\Vert_{L^p(\ell^q)},
 \end{align*}
 as desired.
  
 \hfill

 \subsection{The case when $0<p<q\le \infty$}
 
 In this case, we assume $p>t$ and apply a duality argument with $1<p/t<q/t\le \infty$.
We write
\begin{align*}
 \Big\Vert \Big\{\mathfrak{M}_{\sigma,2^k,y}^{t}f_k \Big\}_{k\in\bbz}\Big\Vert_{L^p(\ell^q)}&=\bigg\Vert   \bigg\{ \Big( \mathfrak{M}_{\sigma,2^k,y}^{t}f_k \Big)^t \bigg\}_{k\in\bbz}     \bigg\Vert_{L^{p/t}(\ell^{q/t})}^{1/t}\\
 &=\sup_{\Vert \{g_k\}_{k\in\bbz}\Vert_{L^{(p/t)'}(\ell^{(q/t)'})}=1}\bigg| \int_{\bbrn} \sum_{k\in\bbz}   \Big( \mathfrak{M}_{\sigma,2^k,y}^{t}f_k(x) \Big)^t g_k(x)       \; dx\bigg|^{1/t}.
 \end{align*}
 Then 
 \begin{align}\label{maximaldualidea}
 &\bigg| \int_{\bbrn} \sum_{k\in\bbz}   \Big( \mathfrak{M}_{\sigma,2^k,y}^{t}f_k(x) \Big)^t g_k(x)       \; dx \bigg|\nonumber\\
 &\le \int_{\bbrn}\int_{\bbrn}\sum_{k\in\bbz} \frac{2^{kn}}{\big(1+2^k|x-z-2^{-k}y| \big)^{\sigma t}}\big| f_k(z)\big|^t \big| g_k(x)\big| \; dx  dz\nonumber\\
 &=\int_{\bbrn} \sum_{k\in\bbz}\big|f_k(z) \big|^t \mathfrak{M}_{\sigma t,2^k,-y}^1g_k(z) \; dz
 \end{align}
 and by H\"older's inequality with $p/t,q/t>1$, the preceding expression is dominated by
 \begin{align*}
& \big\Vert \{f_k\}_{k\in\bbz}\big\Vert_{L^p(\ell^q)}^t\Big\Vert  \Big\{  \mathfrak{M}_{\sigma t,2^k,-y}^1g_k   \Big\}_{k\in\bbz}      \Big\Vert_{L^{(p/t)'}(\ell^{(q/t)'})}\\
&\lesssim \big\Vert \{f_k\}_{k\in\bbz}\big\Vert_{L^p(\ell^q)}^t  \big(\ln{(e+|y|)} \big)^{\frac{1}{(q/t)'}-\frac{1}{(p/t)'}} \big\Vert \{g_k\}_{k\in\bbz}\big\Vert_{L^{(p/t)'}(\ell^{(q/t)'})}     \\
&=\big\Vert \{f_k\}_{k\in\bbz}\big\Vert_{L^p(\ell^q)}^t\big(\ln{(e+|y|)} \big)^{t(1/p-1/q)} 
 \end{align*}
 where the inequality was already verified in Section \ref{0<q<p<inf} because $1\le (q/t)'<(p/t)'<\infty$ and $\sigma t>n$.
 This concludes 
 $$ \Big\Vert \Big\{\mathfrak{M}_{\sigma,2^k,y}^{t}f_k \Big\}_{k\in\bbz}\Big\Vert_{L^p(\ell^q)}\lesssim \big( \ln{(e+|y|)}\big)^{1/p-1/q}\big\Vert \{f_k\}_{k\in\bbz}\big\Vert_{L^p(\ell^q)}.$$
 
\hfill

 \subsection{The case when $p=\infty$ and $0<q< \infty$}
 In this case, we may apply the same arguments as in Section \ref{0<q<p<inf}.
 The left-hand side of \eqref{pinfqfiniteest} is exactly the same as
 $$\bigg\Vert \mathscr{M}_q^{\sharp,(1)}\Big(\Big\{ \mathfrak{M}_{\sigma,2^k,y}^qf_k \Big\}_{k\in\bbz} \Big)\bigg\Vert_{L^{\infty}(\bbrn)}$$
 and in order to estimate this term, we decompose it into 
$$\mathrm{I}_1^{(\infty,q)}+\mathrm{I}_2^{(\infty,q)}+\mathrm{I}_3^{(\infty,q)}$$ 
 as in \eqref{lpsuppind}, where $\mathrm{I}_{j}^{(\infty,q)}$, $j=1,2,3$, are defined by  replacing the $L^p$ norm in \eqref{lpsuppind} by the $L^{\infty}$ norm.
Then analogous to \eqref{II1pest}, \eqref{II2pest}, and \eqref{II3pest}, we have
\begin{equation*}
\mathrm{I}_1^{(\infty,q)}\lesssim \Big\Vert   \mathscr{M}_q^{\sharp,(1)}\big( \{f_k\}_{k\in\bbz}\big)   \Big\Vert_{L^q(\bbrn)}=\sup_{P\in\mathcal{D}}\bigg(\frac{1}{|P|}\int_P \sum_{k:2^k\ell(P)\ge 1} \big| f_k(x)\big|^q\; dx \bigg)^{1/q},
\end{equation*} 
\begin{equation*}
\mathrm{I}_2^{(\infty,q)}\lesssim \big(\ln{(e+|y|)} \big)^{1/q} \big\Vert \{f_k\}_{k\in\bbz} \big\Vert_{L^{\infty}(\ell^{\infty})},
\end{equation*}
 and
\begin{equation*}
\mathrm{I}_3^{(\infty,q)}\lesssim  \big\Vert \{f_k\}_{k\in\bbz} \big\Vert_{L^{\infty}(\ell^{\infty})}.
\end{equation*}
Finally, \eqref{pinfqfiniteest} follows from \eqref{embeddinginfty}.

 \hfill

 \section{Proof of Theorem \ref{dyadmaxtheorem}}\label{proofthm13}
 Without loss of generality, we may assume  $t=1<p<\infty$ as the general case follows from using the dilation
 \begin{equation}\label{dilationargument}
 \big\Vert \mathscr{M}_t^yf\big\Vert_{L^p(\bbrn)}=\big\Vert \mathscr{M}^y\big(|f|^t\big) \big\Vert_{L^{p/t}}^{1/t}
 \end{equation}
 where $p/t>1$. Moreover, we may consider, instead of $\mathscr{M}^yf(x)$,
 $$\wt{\mathscr{M}}^yf(x):=\sup_{Q\in\mathcal{D}:x\in Q}\bigg|\frac{1}{|Q|}\int_{Q(y)}f(z) \; dz \bigg|$$
 as the desired result follows simply from replacing $f$ by $|f|$.
 Therefore, matters reduce to the inequality
 \begin{equation}\label{dyadicreduction}
 \big\Vert \wt{\mathscr{M}}^yf\big\Vert_{L^p(\bbrn)}\lesssim \big( \ln{(e+|y|)}\big)^{1/p}\Vert f\Vert_{L^p(\bbrn)}.
 \end{equation}
 
 For $k\in\bbz$ and $x\in \bbrn$, we denote by $Q^k_x$ the unique dyadic cube of side-length $2^{-k}$ containing $x$ and by $Q^k_x(y)$ its $y$-shifted cube; that is, $Q^k_x(y):=Q^k_x+2^{-k}y$.
 Then we write
 \begin{equation}\label{mmyfx}
 \wt{\mathscr{M}}^yf(x)=\sup_{k\in\bbz}\big| \mathscr{N}_k^yf(x)\big|= \big\Vert \big\{    \mathscr{N}_k^yf(x)     \big\}_{k\in\bbz}\big\Vert_{\ell^{\infty}}
 \end{equation}
 where 
 \begin{equation}\label{scrnkyfdef}
 \mathscr{N}_k^yf(x):=\frac{1}{|Q_x^k|}\int_{Q_x^k(y)} f(z)   \; dz.
 \end{equation}
 Since 
 $$|x-z-2^{-k}y|\lesssim 2^{-k} \q \text{ for }~z\in Q_x^k(y),$$
 we have
 \begin{equation}\label{nkyfmaximalfunction}
\big|  \mathscr{N}_k^yf(x) \big|\lesssim_{\sigma}\int_{\bbrn}  \frac{2^{kn}}{\big(1+2^k|x-z-2^{-k}y| \big)^{\sigma}} \big| f(z)\big|   \; dz=\mathfrak{M}_{\sigma,2^k,y}^1f(x)
 \end{equation}
 for any $\sigma>n$.
Then it follows that
 \begin{align*}
 \big\Vert \wt{\mathscr{M}}^yf\big\Vert_{L^p(\bbrn)}&=\big\Vert \big\{ \mathscr{N}_k^yf\big\}_{k\in\bbz}\big\Vert_{L^p(\ell^{\infty})}=\sup_{\Vert \{g_k\}_{k\in\bbz}\Vert_{L^{p'}(\ell^1)}=1}\bigg| \int_{\bbrn}\sum_{k\in\bbz}\mathscr{N}_k^yf(x)g_k(x) \;dx\bigg|\\
 &\lesssim_{\sigma} \sup_{\Vert \{g_k\}_{k\in\bbz}\Vert_{L^{p'}(\ell^1)}=1} \int_{\bbrn}\sum_{k\in\bbz}\Big( \mathfrak{M}_{\sigma,2^k,y}^1f(x)\Big) \big|g_k(x)\big| \;dx.
 \end{align*}
By the arguments used in \eqref{maximaldualidea}, the integral in the preceding expression would be
$$\int_{\bbrn}\big| f(x)\big| \sum_{k\in\bbz}\mathfrak{M}_{\sigma,2^k,-y}^1g_k(x) \;dx$$
and this is bounded, via H\"older's inequality, by 
 \begin{align*}
 \Vert f\Vert_{L^p(\bbrn)}\Big\Vert \Big\{ \mathfrak{M}_{\sigma,2^k,-y}^1g_k\Big\}_{k\in\bbz}\Big\Vert_{L^{p'}(\ell^1)}.
 \end{align*}
 Consequently,
 \begin{equation*}
 \big\Vert \wt{\mathscr{M}}^yf\big\Vert_{L^p(\bbrn)}\lesssim \Vert f\Vert_{L^p(\bbrn)}\sup_{\Vert \{g_k\}_{k\in\bbz}\Vert_{L^{p'}(\ell^1)}=1}\Big\Vert \Big\{ \mathfrak{M}_{\sigma,2^k,-y}^1g_k\Big\}_{k\in\bbz}\Big\Vert_{L^{p'}(\ell^1)}.
 \end{equation*}
 Now we apply Proposition \ref{maintheorem} to bound the supremum by a constant times
 $$\big( \ln{(e+|y|)}\big)^{1-1/p'}=\big( \ln{(e+|y|)}\big)^{1/p},$$
which completes the proof of \eqref{dyadicreduction}.

 \hfill

 \section{Key lemmas for the proof of Theorem \ref{vectorvalueddyadic}}\label{keylemmathm14}
 
 In this section, we introduce two lemmas, which will be key tools to establish vector-valued estimates in the proof of Theorem \ref{vectorvalueddyadic}.
 \subsection{Vector-valued extension of boundedness for positive linearizable operators} 
 We first recall a vector-valued inequality for linearizable operators, introduced in \cite[Chap.V, Corollary 1.23]{Ga_Ru1985}.
 An operator $T$, acting on $L^p(\bbrn)$, is called linearizable if there exists a linear operator $\mathscr{T}$, acting on $L^p(\bbrn)$, whose values are $\mathscr{B}$-valued functions for some Banach space $\mathscr{B}$ such that
 $$\big| Tf(x)\big|=\big\Vert \mathscr{T}f(x)\big\Vert_{\mathscr{B}}, \qq f\in L^p(\bbrn).$$
Then the following lemma is extremely useful in extending $L^p$ boundedness of linearizable positive operators to a vector-valued space $L^p(\ell^q)$ for $1\le p\le q$.
\begin{lemma}\cite{Ga_Ru1985}\label{exlptolpq}
Let $1\le p<\infty$ and $p\le q\le \infty$. 
Suppose that $T$ is a linearizable operator satisfying
$$\Vert Tf\Vert_{L^p(\bbrn)}\le \mathcal{A}\Vert f\Vert_{L^p(\bbrn)}.$$
If $T$ is positive (in the sense that $\big|Tf(x)\big|\le Tg(x)$ $\mathrm{a.e.}$ holds provided that $|f(x)|\le g(x)$ $\mathrm{a.e.}$), then
 we have
\begin{equation*}
\big\Vert \big\{ Tf_k\big\}_{k\in\bbz}\big\Vert_{L^p(\ell^q)}\lesssim_{p,q}\mathcal{A} \big\Vert \big\{ f_k\big\}_{k\in\bbz}\big\Vert_{L^p(\ell^q)}.
\end{equation*}
\end{lemma}

\hfill

\subsection{Calder\'on-Zygmund theory for vector-valued function spaces}
 For each $j\in\bbz$, $y\in\bbrn$, and $\sigma>n$, let
$$\Lambda_{j,\sigma}^y(x):=\frac{2^{jn}}{(1+|2^jx-y|)^{\sigma}}.$$
  \begin{lemma}\label{keylemma}
 Let $1<p\le q<\infty$, $y\in \bbrn$, and $\sigma>n$.
 Then we have
  \begin{equation}\label{weak11est}
  \bigg\Vert \bigg(\sum_{k\in\bbz}\Big| \sum_{j\in\bbz}\Lambda_{j,\sigma}^y \ast f_{j,k}\Big|^q \bigg)^{1/q}\bigg\Vert_{L^{1,\infty}(\bbrn)}\lesssim  \ln{(e+|y|)} \bigg\Vert \bigg(\sum_{k\in\bbz}   \Big( \sum_{j\in\bbz}\big| f_{j,k}\big|   \Big)^q       \bigg)^{1/q}\bigg\Vert_{L^1(\bbrn)}
  \end{equation}
 and
 \begin{align}\label{strongpq1est}
 &\bigg\Vert \bigg(\sum_{k\in\bbz}\Big| \sum_{j\in\bbz}\Lambda_{j,\sigma}^y\ast f_{j,k}\Big|^q \bigg)^{1/q}\bigg\Vert_{L^p(\bbrn)}\nonumber\\
 &\lesssim \big(\ln{(e+|y|)}\big)^{1-\frac{p-1}{p(q-1)}}\bigg\Vert \bigg(\sum_{k\in\bbz}   \Big( \sum_{j\in\bbz}\big| f_{j,k}\big|   \Big)^q       \bigg)^{1/q}\bigg\Vert_{L^p(\bbrn)}.
 \end{align}
 \end{lemma}
 \begin{proof}
 We first observe that
  $$\big|\Lambda_{j,\sigma}^y\ast f_{j,k}(x)\big|\le \mathfrak{M}_{\sigma,2^j,y}^1f_{j,k}(x).$$
  When $p=q$, it follows from Proposition \ref{maintheorem}  that
\begin{align}\label{strongqqest}
\bigg\Vert \bigg(\sum_{k\in\bbz}\Big| \sum_{j\in\bbz}\Lambda_{j,\sigma}^y\ast f_{j,k}\Big|^q \bigg)^{1/q}\bigg\Vert_{L^q(\bbrn)}&\le  \bigg(\sum_{k\in\bbz}\Big\Vert \Big\{ \mathfrak{M}_{\sigma,2^j,y}^1 f_{j,k}\Big\}_{j\in\bbz}\Big\Vert_{L^q(\ell^1)}^q \bigg)^{1/q}\nonumber\\
&\lesssim \big( \ln{(e+|y|)}\big)^{1-1/q} \bigg(\sum_{k\in\bbz}\big\Vert \big\{  f_{j,k}\big\}_{j\in\bbz}\big\Vert_{L^q(\ell^1)}^q \bigg)^{1/q}\nonumber\\
&= \big( \ln{(e+|y|)}\big)^{1-1/q} \bigg\Vert \bigg(\sum_{k\in\bbz}   \Big( \sum_{j\in\bbz}\big| f_{j,k}\big|   \Big)^q       \bigg)^{1/q}\bigg\Vert_{L^q(\bbrn)},
\end{align}
 as desired.

Moreover, 
 by taking
 \begin{equation}\label{aysupwinrn}
 \mathcal{A}_y:=\sup_{w\in\bbrn}\int_{|x|>2|w|}\sup_{j\in\bbz}\big| \Lambda_{j,\sigma}^y(x-w)-\Lambda_{j,\sigma}^y(x)\big|\; dx
 \end{equation}
 and using the Calder\'on-Zygmund theory for vector-valued function spaces, together with \eqref{strongqqest},
 we have
 \begin{align}\label{weak11eststt}
&\bigg\Vert \bigg(\sum_{k\in\bbz}\Big| \sum_{j\in\bbz}\Lambda_{j,\sigma}^y\ast f_{j,k}\Big|^q \bigg)^{1/q}\bigg\Vert_{L^{1,\infty}(\bbrn)}\nonumber\\
&\lesssim \Big( \big( \ln{(e+|y|)}\big)^{1-1/q}+ \mathcal{A}_y\Big) \bigg\Vert \bigg(\sum_{k\in\bbz}   \Big( \sum_{j\in\bbz}\big| f_{j,k}\big|   \Big)^q       \bigg)^{1/q}\bigg\Vert_{L^1(\bbrn)}.
\end{align}
Such a theory appears in many references such as \cite[Theorem 5.17]{Du0}, \cite[Theorem 3.4 in Chapter V]{Ga_Ru1985}, and \cite[Theorem 4.6.1]{Gr0}, but we contain this in the appendix for the sake of completeness.
Now we claim that
\begin{equation}\label{aayest}
\mathcal{A}_y\lesssim \ln{(e+|y|)}.
\end{equation}
We may actually replace $\sup_{w\in\bbrn}$ in \eqref{aysupwinrn} by $\sup_{w\in\bbrn\setminus \{0\}}$ as the integrand vanishes when $w=0$.
 For a fixed $w\in \bbrn\setminus \{0\}$, 
\begin{align*}
&\int_{|x|>2|w|}\sup_{j\in\bbz}\big| \Lambda_{j,\sigma}^y(x-w)-\Lambda_{j,\sigma}^y(x)\big|\; dx\\
&\le \sum_{j\in\bbz}\int_{|x|>2|w|}\big| \Lambda_{j,\sigma}^y(x-w)-\Lambda_{j,\sigma}^y(x)\big|\; dx\\
&\le \sum_{j:2^j|w|<2^{-2}}\cdots + \sum_{j:2^{-2}\le 2^j|w|<2^2(e+|y|)}\cdots  +    \sum_{j:  2^2(e+|y|)\le 2^j|w|}\cdots\\
&=:\II_{y,w}^1+\II_{y,w}^2+\II_{y,w}^3.
\end{align*} 
  To estimate $\II_{y,w}^1$, we see that
  \begin{align*}
 & \int_{|x|>2|w|}\big| \Lambda_{j,\sigma}^y(x-w)- \Lambda_{j,\sigma}^y(x)\big|\; dx  \\
 &\lesssim 2^j |w|\int_0^1\int_{|x|>2|w|}\frac{2^{jn}}{(1+|2^jx-y-2^jwt|)^{\sigma+1}} \; dx dt\lesssim 2^j|w|
  \end{align*}
 and thus
 $$\II_{y,w}^1\lesssim \sum_{j:2^j|w|<2^{-2}}2^j|w|\lesssim 1.$$
 Moreover, since 
 \begin{equation*}
 \int_{|x|>2|w|}\big|  \Lambda_{j,\sigma}^y(x-w)- \Lambda_{j,\sigma}^y(x)\big|\; dx\le \int_{\bbrn}\big|  \Lambda_{j,\sigma}^y(x-w)\big|\; dx+\int_{\bbrn}\big| \Lambda_{j,\sigma}^y(x)\big|\; dx \lesssim 1,
 \end{equation*}
 we have
 $$\II_{y,w}^2\lesssim    \sum_{j:2^{-2}\le 2^j|w|<2^2(e+|y|)} 1 \lesssim \ln{(e+|y|)}.$$
 For the remaining term $\II_{y,w}^3$, we observe that if $2^2(e+|y|)\le 2^j|w|$ and $2|w|<|x|$, then
 $$|2^j(x-w)+y|, |2^jx+y|\gtrsim 2^j|x|,$$
 and thus,
 \begin{equation*}
\int_{|x|>2|w|} \big| \Lambda_{j,\sigma}^y(x-w)-\Lambda_{j,\sigma}^y(x)\big|\; dx\lesssim_{\sigma}   2^{-j(\sigma-n)}\int_{|x|>2|w|}|x|^{-\sigma}\; dx\sim (2^j|w|)^{-(\sigma-n)}.
 \end{equation*}
 This yields
 $$\II_{y,w}^3\lesssim  \sum_{j:  2^2(e+|y|)\le 2^j|w|} (2^{j}|w|)^{-(\sigma-n)}\lesssim 1,$$
 which completes the proof of  \eqref{aayest}.
Now \eqref{weak11est} follows from \eqref{weak11eststt}. 
 
 Finally, when $1<p<q$, \eqref{strongpq1est}  follows from interpolating \eqref{strongqqest} and \eqref{weak11est} with
 $$1-\frac{p-1}{p(q-1)}=\Big(1-\frac{1}{q}\Big)\theta+1-\theta$$
 where $0< \theta<1$ is the number satisfying
 $\frac{1}{p}=\frac{\theta}{q}+\frac{1-\theta}{1}.$
 \end{proof}

 \hfill

 \section{Proof of Theorem \ref{vectorvalueddyadic}}\label{proofthm14}
 \subsection{The case when $p\le q\le \infty$}
We may assume $1=t<p\le q\le \infty$ as the estimate in the full range $0<t<p\le q\le \infty$ follows from the same dilation argument as in \eqref{dilationargument}. Moreover, it is enough to obtain estimates for  $\wt{\mathscr{M}}^y$, instead of $\mathscr{M}^y$, similar to the reduction in \eqref{dyadicreduction}.
 We appeal to Lemma \ref{exlptolpq}.
For each $y\in\bbrn$ and $k\in\bbz$, 
$\mathcal{N}_k^y$, defined in \eqref{scrnkyfdef}, is a linear operator and then as discussed in \cite[page 481, (1.21)]{Ga_Ru1985},
$\wt{\mathscr{M}}^y$ is linearizable when we take $\mathscr{B}=\ell^{\infty}$.
Moreover, if $|f(x)|\le g(x)$ $\mathrm{a.e.}$,
then 
\begin{align*}
\wt{\mathscr{M}}^yf(x)\le \sup_{k\in\bbz}\frac{1}{|Q_x^k|}\int_{Q_x^k(y)}|f(z)|\; dz\le  \sup_{k\in\bbz}\frac{1}{|Q_x^k|}\int_{Q_x^k(y)}g(z)\; dz=\wt{\mathscr{M}}^yg(x)
\end{align*}
and thus $\wt{\mathscr{M}}$ is a positive operator.
Finally, applying Lemma \ref{exlptolpq} with \eqref{dyadicreduction}, we obtain
 \begin{equation*}
 \big\Vert \big\{ \wt{\mathscr{M}}^yf_k\big\}_{k\in\bbz}\big\Vert_{L^p(\ell^q)}\lesssim \big( \ln{(e+|y|)}\big)^{1/p}\big\Vert \{f_k\}_{k\in\bbz}\big\Vert_{L^p(\ell^q)},
 \end{equation*}
 as desired.
 
\hfill

 \subsection{The case when $q<p<\infty$}
In this case, the following proposition is the main estimate, which is actually weaker inequality (with greater logarithmic exponent) than the desired one \eqref{vecdyashiftest}, but can be refined by a suitable dilation argument. 
\begin{proposition}\label{firstreduction}
Let $1<q<p<\infty$.
Then we have
 \begin{equation}\label{firstreductionest}
 \big\Vert   \big\{ {\mathscr{M}}^yf_k\big\}_{k\in\bbz}  \big\Vert_{L^p(\ell^q)}\lesssim \big( \ln{(e+|y|)}\big)^{1-\frac{q}{p}+\frac{1}{p}}\big\Vert \big\{f_k \big\}_{k\in\bbz}\big\Vert_{L^p(\ell^q)}.
 \end{equation}
\end{proposition}
We postpone the proof of the proposition to the end of this subsection and prove Theorem \ref{vectorvalueddyadic}.
 Suppose that $0<t<q<p<\infty$ and let $\tau_{p,q}$ be any number greater than $1/q$.
 Observing
 $$ \frac{1}{p}+\frac{q}{s}\Big(\frac{1}{q}-\frac{1}{p} \Big) \searrow \frac{1}{q} \q\text{ as }\q {s\nearrow q},$$
 we choose $t<s<q$ such that
 \begin{equation}
 \tau_{p,q}>\frac{1}{p}+\frac{q}{s}\Big(\frac{1}{q}-\frac{1}{p} \Big).
 \end{equation}
 Finally, using Proposition \ref{firstreduction}, we deduce that
 \begin{align*}
  \big\Vert   \big\{ {\mathscr{M}}_t^yf_k\big\}_{k\in\bbz}  \big\Vert_{L^p(\ell^q)}&\le   \big\Vert   \big\{ {\mathscr{M}}_s^yf_k\big\}_{k\in\bbz}  \big\Vert_{L^p(\ell^q)}=   \Big\Vert   \Big\{ {\mathscr{M}}^y\big(|f_k|^s\big)\Big\}_{k\in\bbz}  \Big\Vert_{L^{p/s}(\ell^{q/s})}^{1/s}\\
  &\lesssim_{p,q,s} \bigg( \big( \ln{(e+|y|)}\big)^{1-\frac{q}{p}+\frac{s}{p}}    \Big\Vert \Big\{ |f_k|^s\Big\}_{k\in\bbz}\Big\Vert_{L^{p/s}(\ell^{q/s})}     \bigg)^{1/s}\\
  &= \big(\ln{(e+|y|)} \big)^{\frac{1}{p}+\frac{q}{s}(\frac{1}{q}-\frac{1}{p})} \big\Vert \big\{f_k \big\}_{k\in\bbz}\big\Vert_{L^p(\ell^q)}\\
  &\le \big(\ln{(e+|y|)} \big)^{\tau_{p,q}} \big\Vert \big\{f_k \big\}_{k\in\bbz}\big\Vert_{L^p(\ell^q)},
 \end{align*}
 as desired.

\begin{proof}[Proof of Proposition \ref{firstreduction}]
 Without loss of generality, we may assume that
 $f_k\ge 0$ for all $k\in\bbz$ and then $${\mathscr{M}}^yf_k=\wt{\mathscr{M}}^yf_k.$$
 In this case, in view of \eqref{nkyfmaximalfunction} and \eqref{mmyfx},
 $$\mathscr{N}_j^yf_k(x)\lesssim_{\sigma}\mathfrak{M}_{\sigma,2^j,y}^1f_k(x)=\Lambda_{j,\sigma}^y \ast f_k(x)$$
 and
 $$\wt{\mathscr{M}}^yf_k(x)=\big\Vert \big\{ \mathscr{N}_j^yf_k(x)\big\}_{j\in\bbz}\big\Vert_{\ell^{\infty}}\lesssim \big\Vert \big\{ \Lambda_{j,\sigma}^y \ast f_k(x)\big\}_{j\in\bbz}\big\Vert_{\ell^{\infty}}.$$
 Therefore, \eqref{firstreductionest} can be reduced to the estimate
\begin{equation}\label{secondreduction}
 \Big\Vert   \Big\{ \big\Vert \big\{ \Lambda_{j,\sigma}^y \ast f_k\big\}_{j\in\bbz}\big\Vert_{\ell^{\infty}}\Big\}_{k\in\bbz}  \Big\Vert_{L^p(\ell^q)}\lesssim \big( \ln{(e+|y|)}\big)^{1-\frac{q}{p}+\frac{1}{p}}\big\Vert \big\{f_k \big\}_{k\in\bbz}\big\Vert_{L^p(\ell^q)}.
\end{equation} 
By duality, the left-hand side of \eqref{secondreduction} is equal to
$$\sup_{\Vert \{g_{j,k}\}_{j,k\in\bbz}\Vert_{L^{p'}(\ell^{q'}(\ell^1))}=1}\bigg| \int_{\bbrn} \sum_{k\in\bbz}\sum_{j\in\bbz}   \Lambda_{j,\sigma}^y \ast f_k(x) g_{j,k}(x)\; dx     \bigg|$$
where
$$\big\Vert \big\{g_{j,k}\big\}_{j,k\in\bbz}\big\Vert_{L^{p'}(\ell^{q'}(\ell^1))}:= \Big\Vert   \Big\{ \big\Vert \big\{g_{j,k}(x)\big\}_{j\in\bbz}\big\Vert_{\ell^{1}}\Big\}_{k\in\bbz}  \Big\Vert_{L^{p'}(\ell^{q'})}.$$
By using H\"older's inequality, we have
 \begin{align*}
& \bigg| \int_{\bbrn} \sum_{k\in\bbz}\sum_{j\in\bbz}   \Lambda_{j,\sigma}^y \ast f_k(x) g_{j,k}(x)\; dx     \bigg|\\
&=\bigg| \int_{\bbrn} \sum_{k\in\bbz}\sum_{j\in\bbz}    f_k(x) \Lambda_{j,\sigma}^{-y} \ast g_{j,k}(x)\; dx     \bigg|\\
 &\le \big\Vert \big\{f_k \big\}_{k\in\bbz}\big\Vert_{L^p(\ell^q)} \bigg\Vert \bigg(\sum_{k\in\bbz}\Big( \sum_{j\in\bbz}\Lambda_{j,\sigma}^{-y}\ast |g_{j,k}| \Big)^{q'} \bigg)^{1/q'}\bigg\Vert_{L^{p'}(\bbrn)}
 \end{align*}
 and then Lemma \ref{keylemma} yields that
  the $L^{p'}$ norm in the preceding expression is no greater than a constant multiple of
 $$\big( \ln{(e+|y|)}\big)^{1-\frac{p'-1}{p'(q'-1)}}\big\Vert \big\{g_{j,k}\big\}_{j,k\in\bbz}\big\Vert_{L^{p'}(\ell^{q'}(\ell^1))}=\big( \ln{(e+|y|)}\big)^{1-\frac{q}{p}+\frac{1}{p}}.$$
 This concludes the proof of \eqref{secondreduction}.
 \end{proof}
 
\hfill
 
 \section{Counterexamples}\label{counterexamplesec}
Throughout the section, let $\eta$ be a nonnegative Schwartz function such that $\eta(x)\ge c$ on $\{x\in\bbrn:|x|\le 1/100\}$ for some $c>0$, and $\supp(\wh{\eta})\subset \{\xi\in\bbrn: |\xi|\le 1/1000\}$.
  Without loss of generality, we may assume that $|y|>10e$.

 \subsection{Proof of Theorem \ref{maincorol} : ``if" part}
 Without loss of generality, we may assume $A=1$, which actually affected the results only up to a constant.

 \subsubsection{When $0<p<q\le \infty$}

 Let $$f_k(x):=\begin{cases}
 \eta(x)e^{2\pi i\langle x,2^k e_1\rangle}, & 1\le k\le \ln{(e+|y|)}\\
 0,& otherwise
 \end{cases}$$
 Then we first observe that $f_k\in\mathcal{E}(2^k)$ as
 $$\supp{\wh{f_k}}\subset \{\xi\in\bbrn: 2^k-1/1000\le |\xi|\le 2^k+1/1000\}$$
 for $1\le k\le \ln{(e+|y|)}$.
 First of all,
 \begin{align*}
 \big\Vert \big\{ f_k\big\}_{k\in\bbz}\big\Vert_{L^p(\ell^q)}&=\big\Vert \big\{ \eta \big\}_{1\le k\le \ln{(e+|y|)}}\big\Vert_{L^p(\ell^q)}\\
 &\lesssim \big(\ln{(e+|y|)}\big)^{1/q}\Vert \eta\Vert_{L^p(\bbrn)}\sim \big(\ln{(e+|y|)}\big)^{1/q}.
 \end{align*}
 On the other hand, using \eqref{decreasingshiftn},
\begin{equation}\label{lowerestkey}
\mathfrak{M}_{\sigma,2^k,y}^tf_k(x)\gtrsim \mathfrak{M}_{\sigma,2^k}f_k(x-2^{-k}y)\ge \big|f_k(x-2^{-k}y)\big|=\big| \eta(x-2^{-k}y)\big|
\end{equation}
for $1\le k\le \ln{(e+|y|)}$, and thus
\begin{align*}
 &\Big\Vert \Big\{ \mathfrak{M}_{\sigma,2^k,y}^tf_k\Big\}_{k\in\bbz}\Big\Vert_{L^p(\ell^q)} \\
 &\gtrsim \Big\Vert  \big\{ \eta(\cdot-2^{-k}y)\big\}_{1\le k\le \ln{(e+|y|)}}\Big\Vert_{L^p(\ell^q)}\\
 &\ge \bigg( \sum_{1\le m\le \ln{(e+|y|)}}\int_{|x-2^{-m}y|\le 1/100}\big\Vert \big\{ \eta(x-2^{-k}y)\big\}_{1\le k\le \ln{(e+|y|)}}\big\Vert_{\ell^q}^p \; dx\bigg)^{1/p}\\
 &\ge  \bigg( \sum_{1\le m\le \ln{(e+|y|)}}\int_{|x-2^{-m}y|\le 1/100} \big| \eta(x-2^{-m}y)\big|^p \; dx\bigg)^{1/p}\\
 &\gtrsim \bigg(\sum_{1\le m\le \ln{(e+|y|)}}  1 \bigg)^{1/p}\sim \big(\ln{(e+|y|)}\big)^{1/p}.
\end{align*}
 Therefore, if \eqref{finitepest} holds, then we conclude that
 $$\big( \ln{(e+|y|)}\big)^{1/p}\lesssim \big( \ln{(e+|y|)}\big)^{r+1/q}\q \text{ uniformly in }~ y\in\bbrn.$$
 This yields that
 $$r\ge 1/p-1/q.$$
 
 \subsubsection{When $0<q\le p<\infty$}

 Let 
 \begin{equation}\label{fkxconstruct}
 f_k(x):=\begin{cases}
 \eta(x+2^{-k}y)e^{2\pi i\langle x,2^k e_1\rangle}, & 1\le k\le \ln{(e+|y|)}\\
 0, & otherwise
 \end{cases}.
 \end{equation}
 Then it is easy to verify $f_k\in\mathcal{E}(2^k)$ and
  \begin{equation*}
  \big\Vert \big\{ f_k\big\}_{k\in\bbz}\big\Vert_{L^p(\ell^q)}=  \Big\Vert  \big\{ \eta(\cdot+2^{-k}y)\big\}_{1\le k\le \ln{(e+|y|)}} \Big\Vert_{L^p(\ell^q)}
 \end{equation*}
 The preceding expression is bounded by a constant multiple of 
 $$\bigg\Vert \bigg(\sum_{1\le k\le \ln{(e+|y|)}}\frac{1}{(1+|\cdot+2^{-k}y|)^{2qM}} \bigg)^{1/q}\bigg\Vert_{L^p(\bbrn)}$$
 for sufficiently large $M>n/q$, 
 and this is controlled by
 $$\bigg( \int_{\bbrn}\sum_{1\le k\le \ln{(e+|y|)}}\frac{1}{(1+|x+2^{-k}y|)^{pM}}\; dx\bigg)^{1/p}\sim \big(\ln{(e+|y|)}\big)^{1/p}$$
by using H\"older's inequality if $p/q>1$ and the fact that
 \begin{equation}\label{lowerkeyest1}
 \sum_{1\le k\le \ln{(e+|y|)}}\frac{1}{(1+|x+2^{-k}y|)^L}\lesssim_L 1\q \text{uniformly in ~$x$}
 \end{equation}
 for $L>1$.

Moreover, similar to \eqref{lowerestkey}, we have
\begin{equation}\label{lowerpointest}
\mathfrak{M}_{\sigma,2^k,y}^tf_k(x)\gtrsim \mathfrak{M}_{\sigma,2^k}f_k(x-2^{-k}y)\ge \big| f_k(x-2^{-k}y)\big|=\big| \eta(x)\big|
\end{equation}
for $1\le k\le \ln{(e+|y|)}$.
This yields that
\begin{align*}
\Big\Vert \Big\{ \mathfrak{M}_{\sigma,2^k,y}^tf_k\Big\}_{k\in\bbz}\Big\Vert_{L^p(\ell^q)} &\gtrsim \bigg(\int_{\bbrn}|\eta(x)|^p\bigg( \sum_{1\le k\le \ln{(e+|y|)}}1\bigg)^{p/q} \; dx\bigg)^{1/p}\\
&\sim \Vert \eta\Vert_{L^p(\bbrn)} \big( \ln{(e+|y|)}\big)^{1/q}\sim  \big( \ln{(e+|y|)}\big)^{1/q}.
\end{align*}
 
 Applying these estimates to \eqref{finitepest}, it holds that
 $$\big( \ln{(e+|y|)}\big)^{1/q}\lesssim \big( \ln{(e+|y|)}\big)^{r+1/p}\q \text{ uniformly in }~ y\in\bbrn,$$
 which verifies
 $$r\ge 1/p-1/q.$$

\subsubsection{When $p=q=\infty$}
 In this case, it can be verified with the same construction used just before.
 Indeed, taking \eqref{fkxconstruct}, it is clear that
 $$  \big\Vert \big\{ f_k\big\}_{k\in\bbz}\big\Vert_{L^{\infty}(\ell^{\infty})}=  \Big\Vert  \big\{ \eta(\cdot+2^{-k}y)\big\}_{1\le k\le \ln{(e+|y|)}} \Big\Vert_{L^{\infty}(\ell^{\infty})}=\Vert \eta\Vert_{L^{\infty}(\bbrn)}\sim 1,$$
 and using \eqref{lowerpointest},
 $$ \Big\Vert \Big\{ \mathfrak{M}_{\sigma,2^k,y}^tf_k\Big\}_{k\in\bbz}\Big\Vert_{L^{\infty}(\ell^{\infty})}\gtrsim \Vert \eta\Vert_{L^{\infty}(\bbrn)}\sim 1.$$
 Plugging these estimates into \eqref{finitepest}, we obtain
 $$r\ge 0,$$
 as desired.

 \subsubsection{When $0<q<p=\infty$}
 Let $f_k\in\mathcal{E}(2^k)$ be defined as in \eqref{fkxconstruct}.
 Then \eqref{lowerkeyest1} proves
 \begin{align*}
 \big\Vert \{f_k(x)\}_{k\in\bbz}\big\Vert_{\ell^q}&=\bigg( \sum_{1\le k\le \ln{(e+|y|)}}\big|\eta(x+2^{-k}y) \big|^q\bigg)^{1/q}\\
 &\lesssim_L \bigg( \sum_{1\le k\le \ln{(e+|y|)}}\frac{1}{(1+|x+2^{-k}y|)^{Lq}}\bigg)^{1/q}\lesssim 1
 \end{align*} and thus
 \begin{align*}
 \sup_{P\in\mathcal{D}}\bigg( \frac{1}{|P|}\int_P{\sum_{k:2^k\ell(P)\ge 1}}\big| f_k(x)\big|^q\; dx\bigg)^{1/q}\le \sup_{P\in\mathcal{D}}\bigg( \frac{1}{|P|}\int_P{\big\Vert \{f_k(x)\}_{k\in\bbz}\big\Vert_{\ell^q}^q}\; dx\bigg)^{1/q}\lesssim 1.
 \end{align*}
 Moveover, using \eqref{lowerpointest},
\begin{align*}
 &\sup_{P\in\mathcal{D}}\bigg( \frac{1}{|P|}\int_P \sum_{k:2^k\ell(P)\ge 1}\Big(\mathfrak{M}_{\sigma,2^k,y}^{t}f_k(x) \Big)^q\;dx \bigg)^{1/q} \\
 &\ge \bigg( \int_{[0,1]^n} \sum_{k=0}^{\infty} \Big(\mathfrak{M}_{\sigma,2^k,y}^{t}f_k(x) \Big)^q\; dx\bigg)^{1/q}\\
 &\gtrsim \bigg( \int_{[0,1]^n} |\eta(x)|^q  \Big( \sum_{1\le k\le \ln{(e+|y|)}}1 \Big)   \; dx\bigg)^{1/q}\sim \big(\ln{(e+|y|)}\big)^{1/q}.
 \end{align*}
 
 Finally,  if \eqref{pinfqfiniteesto} holds, then we have
 $$\big(\ln{(e+|y|)}\big)^{1/q}\lesssim \big( \ln{(e+|y|)}\big)^{r},$$
 which completes the validity of \eqref{rcondition}.
 
 \hfill
 
 \subsection{Sharpness of Theorem \ref{phikyastfest}}

 Let $f:=\psi$
 so that
 $$\Vert f\Vert_{H^p(\bbrn)}\sim 1.$$
For each $m\ge 1$,
 $$\phi_m\ast f(x)=\phi_m\ast\psi(x)= \psi(x)$$
 and this implies that
 $$(\phi_m)^y\ast f(x)=\psi(x-2^{-m}y).$$
 Therefore,
 \begin{align*}
 \Big\Vert \sup_{k\in\bbz}\big| (\phi_k)^y\ast f\big| \Big\Vert_{L^p(\ell^{\infty})}&\ge \bigg(\sum_{1\le m\le \ln{(e+|y|)}}\int_{|x-2^{-m}y|\le \frac{1}{100}}  \sup_{k\in\bbz}\big| (\phi_k)^y\ast f\big|^p\; dx\bigg)^{1/p}\\
 &\ge \bigg(\sum_{1\le m\le \ln{(e+|y|)}}\int_{|x-2^{-m}y|\le \frac{1}{100}}  \big| (\phi_m)^y\ast f\big|^p\; dx\bigg)^{1/p}\\
 &=\bigg(\sum_{1\le m\le \ln{(e+|y|)}}\int_{|x-2^{-m}y|\le \frac{1}{100}}  \big| \psi(x-2^{-m}y)\big|^p\; dx\bigg)^{1/p}\\
 &=\Vert \psi\Vert_{L^p(|x|\le \frac{1}{100})}\bigg(\sum_{1\le m\le \ln{(e+|y|)}} 1\bigg)^{1/p}\sim \big( \ln{(e+|y|)}\big)^{1/p}.
 \end{align*}

 \hfill
 
 \subsection{Sharpness of Theorem \ref{shiftedvectorest}}
 
 For each $k\in\bbz$, let $\zeta_k=10k$.
 
  Let $$\wt{\psi}:=2^{-1}\psi(2^{-1}\cdot)+\psi+2\psi(2\cdot) \q\text{ and }\q \wt{\psi}_k:=2^{kn}\wt{\psi}(2^k\cdot)$$
  so that
  \begin{equation}\label{psikidentity}
  \wh{\wt{\psi}_k}(\xi)=1 \q \text{ for }~ 2^{k-1}\le |\xi|\le 2^{k+1}.
  \end{equation}
We observe that for each $k\in\bbz$
 $$(\wt{\psi}_{k})^y(x):=\wt{\psi}_{k}(x-2^{-k}y)=(\psi_{k-1})^{2^{-1}y}(x)+(\psi_k)^y(x)+(\psi_{k+1})^{2y}(x).$$
and thus
\begin{align*}
(\wt{\psi}_{k})^y\ast f(x)=(\psi_{k-1})^{2^{-1}y}\ast f(x)+(\psi_{k})^y\ast f(x)+(\psi_{k+1})^{2y}\ast f(x).
\end{align*}
 Then \eqref{psikyastflplqest} implies that
 \begin{align*}
 &\big\Vert \big\{ (\wt{\psi}_k)^y\ast f \big\}_{k\in\bbz}\big\Vert_{L^p(\ell^q)}\\
 &\lesssim_{p,q}  \big\Vert \big\{ ({\psi}_k)^{2^{-1}y}\ast f \big\}_{k\in\bbz}\big\Vert_{L^p(\ell^q)}+\big\Vert \big\{ ({\psi}_k)^y\ast f \big\}_{k\in\bbz}\big\Vert_{L^p(\ell^q)}+\big\Vert \big\{ ({\psi}_k)^{2y}\ast f \big\}_{k\in\bbz}\big\Vert_{L^p(\ell^q)}\\
 &\lesssim \big(\ln{(e+|y|)}\big)^{r} \big\Vert \big\{ \psi_k\ast f\big\}_{k\in\bbz}\big\Vert_{L^p(\ell^q)}.
 \end{align*}
 Therefore, it suffices to show that
 $r\ge |1/p-1/q|$ is a necessary condition for
 \begin{equation}\label{newpsikyastflplqest}
 \big\Vert \big\{ (\wt{\psi}_k)^y\ast f \big\}_{k\in\bbz}\big\Vert_{L^p(\ell^q)}\lesssim  \big(\ln{(e+|y|)}\big)^{r} \big\Vert \big\{ \psi_k\ast f\big\}_{k\in\bbz}\big\Vert_{L^p(\ell^q)},
 \end{equation} instead of \eqref{psikyastflplqest}.
 
 Similarly, in order to verify the sharpness of \eqref{inftypsikyastflplqest}, we may show that the condition $r\ge 1/q$ should be necessarily required for the following inequality to hold:
 \begin{align}\label{newinftypsikyastflplqest}
&\sup_{P\in\mathcal{D}}\bigg( \frac{1}{|P|}\int_P \sum_{k:2^k\ell(P)\ge 1}\big|(\wt{\psi}_k)^y\ast f(x) \big|^q\; dx\bigg)^{1/q}\nonumber\\
&\lesssim \big( \ln{(e+|y|)}\big)^{r}\sup_{P\in\mathcal{D}}\bigg( \frac{1}{|P|}\int_P \sum_{k:2^k\ell(P)\ge 1}\big|\psi_k \ast f(x) \big|^q\; dx\bigg)^{1/q}.
\end{align}

 \subsubsection{{When $0<p<q\le \infty$}}
 Suppose \eqref{newpsikyastflplqest} holds.

 We define
 $$f(x):=\eta(x)\sum_{0\le k\le \frac{1}{10}\ln{(e+|y|)}}e^{2\pi i\langle x,2^{\zeta_k}e_1\rangle}.$$
 Due to the compact support conditions of $\eta$ and $\psi$, we have
 \begin{equation}\label{psijastfxest}
 \psi_j\ast f(x)=\begin{cases}
 \psi_j\ast\big(\eta \,e^{2\pi i\langle \cdot,2^{\zeta_k}\rangle} \big)(x),& ~\zeta_k-1\le j\le \zeta_k+1, ~0\le k\le \frac{1}{10}\ln{(e+|y|)}\\
 0,& \qq\qq otherwise
 \end{cases}
 \end{equation}
 and for each $0\le k\le \frac{1}{10}\ln{(e+|y|)}$, setting $\sigma>n/p$,
 \begin{equation}\label{psijastfmathfrkest}
 \big| \psi_j\ast f(x)\big|\lesssim \mathfrak{M}_{\sigma,2^{\zeta_k}}\eta(x)\q \text{ for }~\zeta_k-1\le j\le \zeta_k+1
 \end{equation}
 so that
 $$\big\Vert \big\{ \psi_j\ast f(x)\big\}_{j=\zeta_k-1}^{\zeta_k+1}\big\Vert_{\ell^q}\lesssim  \mathfrak{M}_{\sigma,2^{\zeta_k}}\eta(x).$$
 Therefore,
 \begin{align}\label{firstupperbound}
 \big\Vert \big\{ \psi_j\ast f \big\}_{j\in\bbz}\big\Vert_{L^p(\ell^q)}&\lesssim \Big\Vert \Big\{  \mathfrak{M}_{\sigma,2^{\zeta_k}}\eta  \Big\}_{0\le k\le \frac{1}{10}\ln{(e+|y|)}} \Big\Vert_{L^p(\ell^q)}\nonumber\\
 &\lesssim \Vert \eta\Vert_{L^p(\bbrn)}\big\Vert \{1\}_{0\le k\le \frac{1}{10}\ln{(e+|y|)}}\big\Vert_{\ell^q}\sim \big(\ln{(e+|y|)} \big)^{1/q}
 \end{align}
 where the second inequality follows from \eqref{maximal1}.
 
 On the other hand, since
 $$B_m^y:=\Big\{x\in\bbrn: |x-2^{-\zeta_m}y|\le \frac{1}{100}\Big\}, ~ 0\le m\le \frac{1}{10}\ln{(e+|y|)},$$
 are pairwise disjoint sets, we have
 \begin{align*}
 \big\Vert \big\{ (\wt{\psi}_j)^y\ast f\big\}_{j\in\bbz}\big\Vert_{L^p(\ell^q)}&\ge \bigg( \sum_{0\le m\le \frac{1}{10}\ln{(e+|y|)}}\int_{B_m^y}\big\Vert \big\{(\wt{\psi}_j)^y\ast f(x) \big\}_{j\in\bbz}\big\Vert_{\ell^q}^p\; dx\bigg)^{1/p}\nonumber\\
 &\ge \bigg( \sum_{0\le m\le \frac{1}{10}\ln{(e+|y|)}}\int_{B_m^y}\big|\wt{\psi}_{\zeta_m}\ast f(x-2^{-\zeta_m}y) \big|^p\; dx\bigg)^{1/p}\nonumber\\
 &=\bigg( \sum_{0\le m\le \frac{1}{10}\ln{(e+|y|)}}\int_{|x|\le \frac{1}{100}}\big|\wt{\psi}_{\zeta_m}\ast f(x) \big|^p\; dx\bigg)^{1/p}.
 \end{align*}
 
 Using \eqref{psikidentity},
 $$\big| \wt{\psi}_{\zeta_m}\ast f(x)\big|=\Big|\wt{\psi}_{\zeta_m}\ast \big( \eta e^{2\pi i\langle \cdot,2^{\zeta_m}e\rangle}\big)(x)\Big|=\big| \eta(x) e^{2\pi i\langle x,2^{\zeta_m}e\rangle}\big|=\big| \eta(x)\big|$$
 and thus
\begin{align}\label{firstlowerbound}
\big\Vert \big\{ (\wt{\psi}_j)^y\ast f\big\}_{j\in\bbz}\big\Vert_{L^p(\ell^q)}&\gtrsim \bigg( \sum_{0\le m\le \frac{1}{10}\ln{(e+|y|)}}\int_{|x|\le \frac{1}{100}} \big| \eta(x)\big|^p\; dx\bigg)^{1/p}\nonumber\\
&\sim \big( \ln{(e+|y|)}\big)^{1/p}.
\end{align}
 
Taking into account \eqref{newpsikyastflplqest}, \eqref{firstupperbound}, and \eqref{firstlowerbound}, we finally obtain
$$\big(\ln{(e+|y|)}\big)^{1/p}\lesssim \big( \ln{(e+|y|)}\big)^{r+1/q} \q \text{uniformly in }~ y\in\bbrn,$$
which completes the proof of 
$$r\ge 1/p-1/q.$$

 \subsubsection{When $0<q\le p< \infty$}
 Assume that \eqref{newpsikyastflplqest} holds.
 
  Let
 \begin{equation}\label{exfxsample}
 f(x):=\sum_{0\le k\le \frac{1}{10}\ln{(e+|y|)}}\eta(x+2^{-\zeta_k}y)e^{2\pi i\langle x,2^{\zeta_k}e_1\rangle}.
 \end{equation}
 Similar to \eqref{psijastfxest} and \eqref{psijastfmathfrkest},
  \begin{equation*}
 \psi_j\ast f(x)=\begin{cases}
 \psi_j\ast\big(\eta(\cdot+2^{-\zeta_k}y) \,e^{2\pi i\langle \cdot,2^{\zeta_k}e_1\rangle} \big)(x),& ~\zeta_k-1\le j\le \zeta_k+1, ~0\le k\le \frac{1}{10}\ln{(e+|y|)}\\
 0,& \qq\qq otherwise
 \end{cases}
 \end{equation*}
 and for $0\le k\le \frac{1}{10}\ln{(e+|y|)}$, setting $\sigma>n/q$,
 \begin{equation}\label{psiconvolest}
 \big| \psi_j\ast f(x)\big|\lesssim \mathfrak{M}_{\sigma,2^{\zeta_k}}\big( \eta(\cdot+2^{-\zeta_k}y)\big)(x) \q \text{ for }~ \zeta_k-1\le j\le \zeta_k+1.
 \end{equation}
 This yields that
 \begin{align}\label{psijastfjzlplqest}
 \big\Vert \big\{ \psi_j\ast f\big\}_{j\in\bbz}\big\Vert_{L^p(\ell^q)}&\lesssim\Big\Vert \Big\{ \mathfrak{M}_{\sigma,2^{\zeta_k}}\big(\eta(\cdot+2^{-\zeta_k}y) \big)\Big\}_{0\le k\le \frac{1}{10}\ln{(e+|y|)}}\Big\Vert_{L^p(\ell^q)}\nonumber\\
 &\lesssim \Big\Vert \big\{ \eta(\cdot+2^{-\zeta_k}y) \big\}_{0\le k\le \frac{1}{10}\ln{(e+|y|)}}\Big\Vert_{L^p(\ell^q)}.
 \end{align}
 Applying the same argument in \eqref{lowerkeyest1}, the preceding expression is controlled by a constant multiple of 
 \begin{align*}
 & \bigg\Vert \bigg( \sum_{0\le k\le \frac{1}{10}\ln{(e+|y|)}}\frac{1}{(1+|\cdot+2^{-\zeta_k}y|)^{2qM}}\bigg)^{1/q}\bigg\Vert_{L^p(\bbrn)}\\
&\lesssim\bigg(\int_{\bbrn}  \sum_{0\le k\le \frac{1}{10}\ln{(e+|y|)}}\frac{1}{(1+|x+2^{-\zeta_k}y|)^{pM}}   \;dx \bigg)^{1/p}\sim \big( \ln{(e+|y|)}\big)^{1/p}.
 \end{align*} for sufficiently large $M>n/q$.

 Moreover, for $0\le k\le \frac{1}{10}\ln{(e+|y|)}$,
 \begin{equation}\label{wtpsizetakyastf}
 \big| (\wt{\psi}_{\zeta_k})^y\ast f(x)\big|=\big| \wt{\psi}_{\zeta_k}\ast f(x-2^{-\zeta_k}y)\big|=\big| \wt{\psi}_{\zeta_k}\ast \big(\eta\, e^{2\pi i\langle \cdot-2^{-\zeta_k}y,2^{\zeta_k}e_1\rangle} \big)\big|=\big| \eta(x)\big|
 \end{equation}
 where the last equality follows from \eqref{psikidentity}.
 This deduces
 \begin{equation}\label{wtpsijyastfjzpqest}
 \begin{aligned}
 \big\Vert \big\{ (\wt{\psi}_j)^y\ast f\big\}_{j\in\bbz}\big\Vert_{L^p(\ell^q)}&\ge\Big\Vert \big\{ (\wt{\psi}_{\zeta_k})^y\ast f\big\}_{0\le k\le \frac{1}{10}\ln{(e+|y|)}}\Big\Vert_{L^p(\ell^q)}\\
 &=\Vert \eta\Vert_{L^p(\bbrn)} \big\Vert \{1\}_{0\le k\le \frac{1}{10}\ln{(e+|y|)}} \big\Vert_{\ell^q}\sim \big(\ln{(e+|y|)}\big)^{1/q}.
 \end{aligned}
 \end{equation}
 
 Finally, the estimate gives
 $$\big(\ln{(e+|y|)} \big)^{1/q}\lesssim \big(\ln{(e+|y|)} \big)^{r+1/p} \q \text{ uniformly in }~y\in\bbrn,$$
 which proves
 $$r\ge 1/q-1/p,$$
 as desired.

  \subsubsection{When $p=q= \infty$}
 
 Setting \eqref{exfxsample}, the estimates \eqref{psijastfjzlplqest} and \eqref{wtpsijyastfjzpqest} are still valid for $p=q=\infty$. Therefore, we have
 $$ \big\Vert \big\{ \psi_j\ast f\big\}_{j\in\bbz}\big\Vert_{L^{\infty}(\ell^{\infty})}\lesssim \Vert \eta\Vert_{L^{\infty}(\bbrn)}$$
 and
 $$ \big\Vert \big\{ (\wt{\psi}_j)^y\ast f\big\}_{j\in\bbz}\big\Vert_{L^{\infty}(\ell^{\infty})}\gtrsim \Vert \eta\Vert_{L^{\infty}(\bbrn)}\sim 1.$$
Now the assumption \eqref{newpsikyastflplqest} yields
 $$r\ge 0.$$
 
   \subsubsection{When $0<q<p= \infty$}
 Let us assume \eqref{newinftypsikyastflplqest}.
 We take $f$ as in \eqref{exfxsample}.
 Then using \eqref{psiconvolest}, \eqref{maximal2inf}, and \eqref{lowerkeyest1}, we have
 \begin{align*}
 &\sup_{P\in\mathcal{D}}\bigg( \frac{1}{|P|}\int_P \sum_{j:2^j\ell(P)\ge 1}\big|\psi_j \ast f(x) \big|^q\; dx\bigg)^{1/q}\\
 &\le \sup_{P\in\mathcal{D}}\bigg( \frac{1}{|P|}\int_P \sum_{\substack{k:2^{\zeta_k}\ell(P)\ge 1,\\ 0\le k\le \frac{1}{10}\ln{(e+|y|)}}}\Big(    \mathfrak{M}_{\sigma,2^{\zeta_k}}\big( \eta(\cdot+2^{-\zeta_k}y)\big)(x)      \Big)^q\; dx\bigg)^{1/q}\\
 &\lesssim \sup_{P\in\mathcal{D}}\bigg( \frac{1}{|P|}\int_P \sum_{\substack{k:2^{\zeta_k}\ell(P)\ge 1,\\ 0\le k\le \frac{1}{10}\ln{(e+|y|)}}}\big|  \big( \eta(\cdot+2^{-\zeta_k}y)\big)(x)      \big|^q\; dx\bigg)^{1/q}\\
 &\lesssim_L  \sup_{P\in\mathcal{D}}\bigg( \frac{1}{|P|}\int_P \sum_{0\le k\le \frac{1}{10}\ln{(e+|y|)}}\frac{1}{(1+|x+2^{-\zeta_k}y|)^{Lq}}\; dx\bigg)^{1/q}\lesssim 1
 \end{align*}
 for $L>1/q$.
 
 In addition, the left-hand side of \eqref{newinftypsikyastflplqest} is no greater than
 \begin{equation*}
\bigg( \int_{[0,1]^n} \sum_{0\le k\le \frac{1}{10}\ln{(e+|y|)}} \big| (\wt{\psi}_{\zeta_k})^y\ast f(x)\big|^q\; dx\bigg)^{1/q}
 \end{equation*}
 and this is, in view of \eqref{wtpsizetakyastf}, equal to
  \begin{equation*}
\bigg( \int_{[0,1]^n} \big|\eta(x)\big|^q\; dx\bigg)^{1/q}\bigg( \sum_{0\le k\le \frac{1}{10}\ln{(e+|y|)}} 1\bigg)^{1/q}\sim \big( \ln{(e+|y|)}\big)^{1/q}.
 \end{equation*}
 
 In conclusion,
 $$\big(\ln{(e+|y|)} \big)^{1/q}\lesssim \big(\ln{(e+|y|)}\big)^{r},\q \text{ uniformly in }~y\in\bbrn,$$
 which proves $r\ge 1/q$.
 
 This completes the proof.

\hfill
 
 \appendix
 \section{Proof of \eqref{weak11eststt}}
 
 Let $1<q<\infty$ and
 suppose that 
 $$\big\Vert  \big\{  f_{j,k}(x)\big\}_{j,k\in\bbz}\big\Vert_{\ell^q(\ell^1)}:=\bigg( \sum_{k\in\bbz}\Big(\sum_{j\in\bbz}\big|f_{j,k}(x)\big| \Big)^{q}\bigg)^{1/q}$$
belongs to $L^1(\bbrn)$. Without loss of generality, we may assume that
$$\Big\Vert \big\Vert  \big\{  f_{j,k}(\cdot)\big\}_{j,k\in\bbz}\big\Vert_{\ell^q(\ell^1)}\Big\Vert_{L^1(\bbrn)}=1.$$
 We fix $\alpha>0$. 
Then it suffices to show that
\begin{align}\label{weakalphaform}
&\Big|\Big\{x\in\bbrn:      \big\Vert  \big\{ \Lambda_{j,\sigma}^y\ast f_{j,k}(x)\big\}_{j,k\in\bbz}\big\Vert_{\ell^q(\ell^1)}>\alpha         \Big\} \Big|\lesssim \frac{1}{\alpha} \Big(\big( \ln{(e+|y|)}\big)^{1-1/q}+\mathcal{A}_y \Big).
\end{align}
To achieve this, let $\gamma>0$ be arbitrary and
let $\mathfrak{G}_{\alpha, \gamma}=\{Q_l\}_{l\in\bbn}$ be the collection of maximal dyadic cubes such that
$$\frac{1}{|Q_l|}\int_{Q_l}  \big\Vert  \big\{ f_{j,k}(x)\big\}_{j,k\in\bbz}\big\Vert_{\ell^q(\ell^1)}   \; dx >\gamma\alpha.$$ 
 Then by the maximality, 
 \begin{equation*}
 \frac{1}{|Q_l|}\int_{Q_l}  \big\Vert  \big\{ f_{j,k}(x)\big\}_{j,k\in\bbz}\big\Vert_{\ell^q(\ell^1)}   \; dx\le 2^n \gamma \alpha \q \text{ for each }~l\in\bbn.
 \end{equation*}
 For each $j,k\in\bbz$, we define the good function $g_{j,k}(x)$ and the bad function $b_{j,k}(x)$ as
 $$b_{j,k}^l(x):=\bigg( f_{j,k}(x)-\frac{1}{|Q_l|}\int_{Q_l} f_{j,k}(z)\; dz\bigg)\chi_{Q_l}(x),$$
 $$b_{j,k}(x):=\sum_{l\in\bbn}b_{j,k}^l(x),$$
 and
 $$g_{j,k}(x):=f_{j,k}(x)-b_{j,k}(x)$$
 so that
 $$f_{j,k}(x)=g_{j,k}(x)+b_{j,k}(x)=g_{j,k}(x)+\sum_{l\in\bbn}b_{j,k}^l(x).$$
 Then we can easily check that 
 \begin{equation}\label{czproperty1}
 \sum_{l\in\bbn}|Q_l|\le \frac{1}{\gamma\alpha},
 \end{equation}
 \begin{equation}\label{czproperty2}
 \supp(b_{j,k}^l)\subset Q_l, \q \int_{Q_l} b_{j,k}^l(x)\; dx=0,\q \Big\Vert \big\Vert  \big\{ b_{j,k}^l(\cdot)\big\}_{j,k\in\bbz}\big\Vert_{\ell^q(\ell^1)} \Big\Vert_{L^1(\bbrn)}\lesssim \gamma\alpha |Q_l|,
 \end{equation}
and
\begin{equation}\label{czproperty3}
\Big\Vert \big\Vert  \big\{ g_{j,k}(\cdot)\big\}_{j,k\in\bbz}\big\Vert_{\ell^q(\ell^1)}   \Big\Vert_{L^1(\bbrn)}\le 1, \q  \Big\Vert \big\Vert  \big\{ g_{j,k}(\cdot)\big\}_{j,k\in\bbz}\big\Vert_{\ell^q(\ell^1)} \Big\Vert_{L^{\infty}(\bbrn)}\lesssim \gamma \alpha.
\end{equation}
 
We first bound the left-hand side of \eqref{weakalphaform} by the sum of 
\begin{align*}
\mathcal{B}_y^{\mathrm{good}}(\alpha)&:=\Big|\Big\{x\in\bbrn:      \big\Vert  \big\{\Lambda_{j,\sigma}^y\ast g_{j,k}(x)\big\}_{j,k\in\bbz}\big\Vert_{\ell^q(\ell^1)}>\frac{\alpha}{2}         \Big\} \Big|,\\
\mathcal{B}_y^{\mathrm{bad}}(\alpha)&:=\Big|\Big\{x\in\bbrn:      \big\Vert  \big\{\Lambda_{j,\sigma}^y\ast b_{j,k}(x)\big\}_{j,k\in\bbz}\big\Vert_{\ell^q(\ell^1)}>\frac{\alpha}{2}         \Big\} \Big|.
\end{align*}
Using Chebyshev's inequality and the strong boundedness \eqref{strongqqest},
 \begin{align*}
 \mathcal{B}_y^{\mathrm{good}}(\alpha)&\lesssim \frac{1}{\alpha^{q}} \Big\Vert \big\Vert  \big\{ \Lambda_{j,\sigma}^y\ast g_{j,k}(\cdot)\big\}_{j,k\in\bbz}\big\Vert_{\ell^q(\ell^1)}\Big\Vert_{L^q(\bbrn)}^q\\
 &\lesssim \frac{1}{\alpha^q}\big(\ln (e+|y|) \big)^{q-1}\Big\Vert \big\Vert  \big\{  g_{j,k}(\cdot)\big\}_{j,k\in\bbz}\big\Vert_{\ell^q(\ell^1)}\Big\Vert_{L^q(\bbrn)}^q.
 \end{align*}
Due to \eqref{czproperty3}, we have
 \begin{align*}
&\Big\Vert \big\Vert  \big\{  g_{j,k}(\cdot)\big\}_{j,k\in\bbz}\big\Vert_{\ell^q(\ell^1)}\Big\Vert_{L^q(\bbrn)}^q\\
&\le \Big\Vert \big\Vert  \big\{  g_{j,k}(\cdot)\big\}_{j,k\in\bbz}\big\Vert_{\ell^q(\ell^1)}\Big\Vert_{L^{\infty}(\bbrn)}^{q-1} \Big\Vert \big\Vert  \big\{  g_{j,k}(\cdot)\big\}_{j,k\in\bbz}\big\Vert_{\ell^q(\ell^1)}\Big\Vert_{L^1(\bbrn)}\\
&\lesssim \gamma^{q-1}\alpha^{q-1},
\end{align*} 
and thus
\begin{equation}\label{bbygoodest}
\mathcal{B}_y^{\mathrm{good}}(\alpha)\lesssim \frac{1}{\alpha} \gamma^{q-1}  \big(\ln (e+|y|) \big)^{q-1}.
\end{equation}
On the other hand,
$$\mathcal{B}_y^{\mathrm{bad}}(\alpha)\le \bigg| \bigcup_{l\in\bbn}Q_l^*\bigg|+ \bigg|\bigg\{x\in \bigg(\bigcup_{l\in\bbn}Q_l^*\bigg)^c:      \big\Vert  \big\{\Lambda_{j,\sigma}^y\ast b_{j,k}(x)\big\}_{j,k\in\bbz}\big\Vert_{\ell^q(\ell^1)}>\frac{\alpha}{2}         \bigg\} \bigg|$$
where $Q_l^*$ is the concentric dilation of $Q_l$ whose side-length is $10\sqrt{n}\ell(Q_l)$,
and
$$\Big| \bigcup_{l\in\bbn}Q_l^*\Big|\lesssim \sum_{l\in\bbn}|Q_l|\le \frac{1}{\gamma \alpha}$$
due to \eqref{czproperty1}.
Now let us estimate
$$\mathcal{B}_{y,\mathrm{out}}^{\mathrm{bad}}(\alpha):= \bigg|\bigg\{x\in \bigg(\bigcup_{l\in\bbn}Q_l^*\bigg)^c:      \big\Vert  \big\{\Lambda_{j,\sigma}^y\ast b_{j,k}(x)\big\}_{j,k\in\bbz}\big\Vert_{\ell^q(\ell^1)}>\frac{\alpha}{2}         \bigg\} \bigg|.$$
Using Chebyshev's inequality and the triangle inequailty
$$   \big\Vert  \big\{\Lambda_{j,\sigma}^y\ast b_{j,k}(x)\big\}_{j,k\in\bbz}\big\Vert_{\ell^q(\ell^1)}\le \sum_{l\in\bbn}    \big\Vert  \big\{\Lambda_{j,\sigma}^y\ast b_{j,k}^l(x)\big\}_{j,k\in\bbz}\big\Vert_{\ell^q(\ell^1)},$$
we obtain
\begin{align*}
\mathcal{B}_{y,\mathrm{out}}^{\mathrm{bad}}(\alpha)&\lesssim \frac{1}{\alpha}\int_{(\cup_{l\in\bbn}Q_l^*)^c}      \big\Vert  \big\{\Lambda_{j,\sigma}^y\ast b_{j,k}(x)\big\}_{j,k\in\bbz}\big\Vert_{\ell^q(\ell^1)}      \; dx\\
&\le \frac{1}{\alpha}\sum_{l\in\bbn}\int_{(Q_l^*)^c}   \big\Vert  \big\{\Lambda_{j,\sigma}^y\ast b_{j,k}^l(x)\big\}_{j,k\in\bbz}\big\Vert_{\ell^q(\ell^1)}\; dx.
\end{align*}
Using the vanishing moment condition of $b_{j,k}^l$ in \eqref{czproperty2},
\begin{align*}
\Lambda_{j,\sigma}^y\ast b_{j,k}^l(x)=\int_{Q_l}\big( \Lambda_{j,\sigma}^y(x-z)-\Lambda_{j,\sigma}^y(x-c_{Q_l})\big)b_{j,k}^l(z)\; dz
\end{align*}
where $c_{Q_l}$ is the center of the dyadic cube $Q_l$, and thus
\begin{align*}
&\int_{(Q_l^*)^c}   \big\Vert  \big\{\Lambda_{j,\sigma}^y\ast b_{j,k}^l(x)\big\}_{j,k\in\bbz}\big\Vert_{\ell^q(\ell^1)}\; dx\\
&\le \int_{(Q_l^*)^c}\int_{Q_l}\bigg(\sum_{k\in\bbz}\Big(\sum_{j\in\bbz}\big|   \Lambda_{j,\sigma}^y(x-z)-\Lambda_{j,\sigma}^y(x-c_{Q_l})    \big| \big|b_{j,k}^l(z)\big| \Big)^q \bigg)^{1/q}\; dz dx\\
&\le \int_{(Q_l^*)^c}\int_{Q_l}\sup_{j\in\bbz} \big|  \Lambda_{j,\sigma}^y(x-z)-\Lambda_{j,\sigma}^y(x-c_{Q_l}) \big|  \big\Vert  \big\{ b_{j,k}^l(z)\big\}_{j,k\in\bbz}\big\Vert_{\ell^q(\ell^1)}\; dz dx\\
&=\int_{Q_l} \big\Vert  \big\{ b_{j,k}^l(z)\big\}_{j,k\in\bbz}\big\Vert_{\ell^q(\ell^1)} \bigg(  \int_{(Q_l^*)^c}   \sup_{j\in\bbz} \big|  \Lambda_{j,\sigma}^y(x-z)-\Lambda_{j,\sigma}^y(x-c_{Q_l}) \big|        \; dx  \bigg)\; dz\\
&\le \sup_{z\in Q_l}\int_{-c_{Q_l}+(Q_l^*)^c}     \sup_{j\in\bbz} \big| \Lambda_{j,\sigma}^y(x-(z-c_{Q_l}))-\Lambda_{j,\sigma}^y(x) \big|     \; dx\; \Big\Vert \big\Vert  \big\{ b_{j,k}^l(\cdot)\big\}_{j,k\in\bbz}\big\Vert_{\ell^q(\ell^1)} \Big\Vert_{L^1(\bbrn)}.
\end{align*}
If $z\in Q_l$ and $x\in -c_{Q_l}+(Q_l^*)^c$, then
$$|x|>2|z-c_{Q_l}|$$
and this implies that
\begin{align*}
&\sup_{z\in Q_l}\int_{-c_{Q_l}+(Q_l^*)^c}     \sup_{j\in\bbz} \big|  \Lambda_{j,\sigma}^y(x-(z-c_{Q_l}))-\Lambda_{j,\sigma}^y(x) \big|     \; dx\le \mathcal{A}_y.
\end{align*}
Using the last one in \eqref{czproperty2}, we obtain
\begin{align*}
\int_{(Q_l^*)^c}   \big\Vert  \big\{\Lambda_{j,\sigma}^y\ast b_{j,k}^l(x)\big\}_{j,k\in\bbz}\big\Vert_{\ell^q(\ell^1)}\; dx \lesssim \gamma \alpha \mathcal{A}_y |Q_l| 
\end{align*}
and thus
$$\mathcal{B}_{y,\mathrm{out}}^{\mathrm{bad}}(\alpha)\lesssim \gamma \mathcal{A}_y \sum_{l\in\bbn}|Q_l|\le \frac{1}{\alpha}\mathcal{A}_y$$
where \eqref{czproperty1} is applied in the last estimate.
Finally,
$$\mathcal{B}_{y}^{\mathrm{bad}}(\alpha)\lesssim \frac{1}{\alpha}\Big(   \frac{1}{\gamma}  +\mathcal{A}_y\Big).$$
Combining this with \eqref{bbygoodest}, the left-hand side of \eqref{weakalphaform} is bounded by a constant times
$$\frac{1}{\alpha}\Big( \gamma^{q-1}\big( \ln{(e+|y|)}\big)^{q-1}   +\frac{1}{\gamma}+\mathcal{A}_y    \Big).$$
Now \eqref{weakalphaform} follows from taking 
$$\gamma=\big( \ln{(e+|y|)}\big)^{-(1-1/q)}.$$

\hfill

\medskip
\noindent {\bf Acknowledgment:} Part of this work was carried out during my research stay at Osaka University. I would like to thank Naohito Tomita for his invitation, hospitality, and very helpful discussion during the stay. I would also like to thank Akihiko Miyachi for interesting discussion at the Research Institute for Mathematical Sciences (RIMS), Kyoto University.



\end{document}